\documentclass{article}
\usepackage[utf8]{inputenc}
\usepackage{hyphenat}

\usepackage[english]{babel}
\usepackage{amsmath,dsfont}
\usepackage[svgnames]{xcolor}
\usepackage{amsfonts,enumitem}
\usepackage{amssymb}
\usepackage{float}
\usepackage{hyperref} 
\usepackage{textcomp} 
\usepackage[backend=biber, sortcites=true]{biblatex}
\usepackage{csquotes}
\usepackage{caption, subfig} 
\usepackage{float, lipsum}
\usepackage{tabularx} 
\usepackage{graphicx,tikz,}
\usepackage{todonotes}
\usepackage{authblk}
\usepackage{siunitx}[table-number-alignment = center]

\usepackage[skins]{tcolorbox}

\newcommand{\tip}{_\mathrm{tip}}

\usepackage{amsmath,amsfonts,bm}

\def\1{\bm{1}}

\def\eps{{\epsilon}}

\def\vtheta{{\bm{\theta}}}

\def\vx{{\bm{x}}}

\DeclareMathAlphabet{\mathsfit}{\encodingdefault}{\sfdefault}{m}{sl}
\SetMathAlphabet{\mathsfit}{bold}{\encodingdefault}{\sfdefault}{bx}{n}

\DeclareMathOperator*{\argmin}{arg\,min}

\DeclareMathOperator{\Id}{I_d}

\newcommand{\EAC}{E_\mathrm{AC}}
\newcommand{\Ed}{E_\mathrm{d}}

\providecommand{\keywords}[1]{\textbf{\textbf{Keywords. }} #1}

\makeatletter
\newcommand{\customlabel}[2]{%
   \protected@write \@auxout {}{\string \newlabel {#1}{{#2}{\thepage}{#2}{#1}{}} }%
   \hypertarget{#1}{#2}
}
\makeatother

\definecolor{babg}{RGB}{220,255,220}   
\definecolor{batext}{RGB}{0,120,0}     

\usepackage[final, commandnameprefix=ifneeded]{changes}
\definechangesauthor[name=CK,color=blue]{CK}
\definechangesauthor[name=MZ,color=orange]{MZ}

\title{
DeepRitzSplit: A Variational Neural Operator Framework for Gradient Flows with Application to Dendritic Growth
}
\author[1]{Beno\^{i}t Appolaire}
\author[2]{Ludovick Gagnon}
\author[1,2]{Chih-Kang Huang}
\author[1]{Miha Zalo\v{z}nik}
\affil[1]{Universit\'{e} de Lorraine, CNRS, IJL, F-54000 Nancy, France}
\affil[2]{Universit\'{e} de Lorraine, CNRS, Inria, IECL, F-54000 Nancy, France}
\date{\today}

\begin{document}
\maketitle

\begin{abstract}
{
Phase-field models are widely used to simulate microstructure evolution during phase transformations, but their numerical solution often remains computationally demanding due to the nonlinear nature of the governing equations and the multiscale character of the physical phenomena. Neural operator surrogates provide a promising alternative, yet existing approaches are typically data-driven, requiring large training datasets and often exhibiting limited robustness when extrapolating to unseen configurations.

In this work, we introduce DeepRitzSplit, an unsupervised neural operator framework for gradient-flow systems based on a variational formulation of energy-stable splitting schemes. The proposed approach learns the time evolution semigroup associated with a convex-concave decomposition of the free-energy functional by minimizing, at each time step, the variational problem underlying the numerical scheme. By construction, the method is designed to inherit the dissipative structure of the underlying gradient flow while avoiding the need for trajectory datasets.

We further introduce a Reaction-Diffusion Neural Operator (RDNO) architecture motivated by the structure of semi-implicit splitting schemes, with separate reaction and diffusion components to improve training efficiency and interpretability.

The approach is first assessed on the Allen-Cahn equation as a controlled benchmark, and subsequently applied to a phase-field model of anisotropic dendritic growth during solidification. In the Allen-Cahn benchmark, DeepRitzSplit shows better generalization than the supervised baseline on the out-of-distribution configurations tested, using 180 training configurations. In the dendritic growth problem, the method captures the overall phase morphology and dendrite tip velocity well, while the tip radius is substantially more sensitive to the training strategy. A subsequent refinement based on the residual of the splitting scheme markedly improves this second-order morphological quantity. The trained operators also retain good predictive capabilities on configurations substantially different from the training set.

These results demonstrate the potential of variationally trained neural operators as physically informed surrogates for complex phase-field models and motivate their use for accelerated simulations of microstructure evolution.
}
\end{abstract}
\keywords{Phase-field modeling; Dendritic Growth; Deep Learning}


\section{Introduction}
\label{sec:introduction}


Phase transformations play an important role in materials processing because they determine the structure and, consequently, the properties of the material. For instance, during additive manufacturing of metal parts, the molten metal first transforms into a polycrystalline solid structure, followed by further phase transformations in the solid state. The final properties of the part depend on the nature, size, and morphology of the microstructure resulting from the different transformations. Detailed modelling of such transformations most often relies on phase-field methods~\cite{Fix1983,Tourret2022}, which represent the phases through one or several continuous fields, or order parameters. The evolution of these fields is typically governed by the gradient flow of a thermodynamic functional~\cite{Plapp2011} and can be coupled to additional physical fields such as temperature, chemical composition, stress, or fluid flow.
Unfortunately, the numerical solution of phase-field models remains computationally demanding. Their governing equations are generally nonlinear, involve multiple spatial and temporal scales, and require fine discretizations to accurately resolve diffuse interfaces. As a consequence, large-scale simulations, uncertainty quantification studies, inverse problems, and process optimization remain challenging in many practical situations~\cite{Takaki2023}. In this context, machine-learning-based surrogate models have recently emerged as a promising route for accelerating phase-field simulations~\cite{choi2024accelerating}.

\paragraph{State of the Art}

Neural operators have attracted considerable attention in recent years as a framework for approximating solution operators of partial differential equations. Unlike conventional neural networks, which learn mappings between finite-dimensional vectors, neural operators learn mappings between function spaces. This property allows them to generalize across different initial conditions, discretizations, and boundary conditions. Most existing neural operators rely on supervised data-driven learning strategies, including DeepONets, Fourier Neural Operators, and related architectures~\cite{bhattacharya2021model,anandkumar2020neural,lu2019deeponet,Li2020FourierNO,patel2021physics}. A comprehensive benchmark of operator-learning architectures for canonical PDEs is provided by Raoni'c et al.~\cite{raonic2023convolutional}.
\\
While these approaches have demonstrated impressive predictive capabilities, they generally require large training datasets generated by numerical simulations. Furthermore, their performance may deteriorate when extrapolating outside the training distribution or when long-time integration is considered. Oommen et al.\cite{oommen2024rethinking} recently highlighted these limitations and proposed hybrid approaches combining neural and classical numerical solvers. In the context of phase-field modelling, Tseng et al.\cite{tseng2023deep} successfully applied neural operators to dendritic growth simulations, but relied on a fully supervised training procedure.
\\
At the same time, \emph{physics-informed} approaches have been developed to reduce the dependence on large training datasets. Physics-Informed Neural Networks (PINNs), introduced by Raissi et al.\cite{raissi2019physics}, incorporate the governing equations directly into the loss function and can therefore be trained without reference solutions. Reviews of recent developments can be found in\cite{cuomo2022scientific}. Applications to phase-field models, including the Allen-Cahn equation, have been reported in~\cite{wight2020solving,mattey2022novel}. However, PINNs typically require retraining for new initial or boundary conditions, which limits their applicability as general-purpose surrogate models.
\\
Physics-Informed Neural Operators (PINOs) aim to combine the strengths of operator learning and physics-informed training. Rather than relying solely on simulation data, they exploit information contained in the governing equations while preserving the ability to generalize across families of solutions. Recent examples include the finite-element-based PINO framework proposed by Yamazaki et al.\cite{yamazaki2025finite} and the finite-difference PINO developed by Geng et al.\cite{geng2024deep} for the Allen-Cahn equation. More closely related to the present work, Li et al.\cite{li2023phase} introduced a variational PINO for the Allen-Cahn and Cahn-Hilliard equations based on the Deep Ritz method\cite{yu2018deep} and a DeepONet architecture. Nevertheless, most existing PINO studies remain focused on relatively simple benchmark equations, while strongly nonlinear phase-field models remain comparatively unexplored.
\\
A common feature of many phase-field models is that they can be formulated as gradient flows of a free-energy functional. Their dynamics are therefore governed not only by the underlying PDEs but also by the dissipation of free energy. In numerical simulations, this structure is often preserved through energy-stable discretization techniques such as convex-concave splitting schemes. Preserving this dissipative structure is essential for obtaining physically meaningful long-time predictions. However, it is only rarely incorporated explicitly into the training procedure of neural operators.

\paragraph{Main Contributions}

In this work, we introduce \emph{DeepRitzSplit}, an unsupervised framework for training neural operators on gradient-flow systems. The proposed approach is based on the variational formulation of energy-stable splitting schemes and is inspired by the Deep Ritz method~\cite{yu2018deep}. Rather than learning from trajectories generated by numerical simulations, the neural operator is trained by minimizing the optimization problem associated with a convex-concave decomposition of the free-energy functional. In this way, the dissipative structure of the underlying gradient flow is directly embedded into the learning procedure.
\\
Inspired by the work of Bretin et al.~\cite{bretin2022learning}, we also introduce a Reaction-Diffusion Neural Operator (RDNO) architecture tailored to semi-implicit phase-field discretizations. The reaction and diffusion contributions are handled by distinct operators, mirroring the decomposition used in the numerical scheme.
\\
The proposed methodology is assessed on two representative phase-field models of increasing complexity. We first consider the Allen-Cahn equation~\cite{allen1979microscopic} as a canonical gradient-flow benchmark. We then investigate the quantitative phase-field model of anisotropic dendritic growth introduced by Karma and Rappel~\cite{karma1998quantitative}, whose mathematical well-posedness was later established in~\cite{miranville2015anisotropic}. Particular attention is paid to out-of-distribution generalization and to the preservation of physically relevant quantities such as phase morphology and dendrite growth characteristics.

The remainder of the paper is organized as follows. Section~\ref{sec:splitting} introduces the DeepRitzSplit \replaced[id=MZ]{neural operator framework.}{framework and the associated neural operator architectures.} Section~\ref{sec:AC} presents the application to the Allen-Cahn equation. Section~\ref{sec:dendritic_growth} considers anisotropic dendritic growth. Finally, the last section summarizes the main findings and discusses perspectives for future developments.

\section{DeepRitzSplit: a variational neural operator framework for gradient flows}
\label{sec:splitting}

Throughout the paper, we use \emph{DeepRitzSplit} for the variational training framework, \emph{splitting scheme} and \emph{SAV scheme} for the numerical schemes, and \emph{RDNO}, \emph{PR-RDNO}, \emph{UNet}, and \emph{FNO} for neural-operator architectures.

Many phase-field models can be formulated as gradient flows and are commonly solved using energy-stable splitting schemes. The central idea of the present work is to exploit the variational formulation of such schemes for training neural operators.
To this end, we first recall in \S\ref{subsec:stable_scheme} a class of semi-implicit splitting schemes introduced by Zaitzeff et al.~\cite{zaitzeff2021high} for general gradient flows. These schemes preserve the dissipative structure of the underlying problem through an associated minimization principle. We then introduce in \S\ref{subsec:DeepRitzSplit} \emph{DeepRitzSplit}, an unsupervised framework that uses this variational formulation directly as a training objective for operator learning.
The problem-specific neural-operator architectures used with this framework are introduced in the corresponding application sections.

%

\subsection{Unconditionally stable semi-implicit scheme for gradient flows}   
\label{subsec:stable_scheme}

Consider the following evolution equation as $L^2$-gradient flow for a given functional, $E : V \to \mathbb{R}$,
\begin{equation}
   u_t = - \nabla_{L^2} E (u)
   \quad\text{on}\quad\Omega \times (0, T)
\label{eq:gradient_flow}
\end{equation} 
where $\Omega$ is a domain of $\mathbb{R}^2$ and $V$ is a Hilbert subspace of  $L^2(\Omega)$.

Solutions to the gradient flow defined by Eq.~\eqref{eq:gradient_flow} have the energy dissipation property: 
\begin{equation}
   \frac{d}{dt} E(u) = \left \langle \nabla_{L^2} E(u), \frac{\partial u}{\partial t}\right\rangle_{L^2}=  - \| \nabla_{L^2} E(u) \|^2_2 \leq 0.
\label{eq:energy_dissipation}
\end{equation}

Zaitzeff et al. \cite{zaitzeff2021high} proposed a class of semi-implicit schemes
for gradient flows \eqref{eq:gradient_flow} that preserve energy stability \eqref{eq:energy_dissipation}.
The energy $E$ is first split into
\[ 
   E(u) = E_1(u) + E_2(u) \quad \forall u \in V.
\]
In the time discretization scheme we then handle $E_1$ implicitly and $E_2$ explicitly: 
\begin{equation}
   \frac{u_{n+1} - u_n}{\mathrm{dt}} = -\nabla_{L^2} E_1(u_{n+1}) - \nabla_{L^2} E_2(u_n)
\label{eq:splitting_scheme}
\end{equation}
with the given time step size $ \mathrm{dt} > 0 $.
The above scheme corresponds to the optimization problem 
\begin{equation}
\begin{aligned}
   u_{n+1} 
   := \argmin_{u \in V}&\bigg(\frac{1}{2 \mathrm{dt}} \| u - u_n \|_2^2 + E_1(u) 
   \\ 
   &\quad + E_2(u_n) + \langle \nabla_{L^2} E_2(u_n), u- u_n \rangle_{L^2}   \bigg).
   \end{aligned}
\label{eq:minimizing_scheme}
\end{equation}

Zaitzeff et al. showed that for $\mathrm{dt} \leq \frac{1}{\Lambda}$, where 
\begin{equation} 
   \Lambda  = \max \left \{ 0 , \max_{u \in V, \|v\| = 1}  D^2E_2(u)(v, v) \right \},
\end{equation}
the energy dissipation property holds for the numerical scheme \eqref{eq:splitting_scheme}: 
\begin{equation}
   E(u_{n+1}) \leq E(u_n).
\label{eq:energy_stability}
\end{equation}
Such a numerical scheme is said to be \emph{energetically stable}.

In particular, in the case where $E_2$ is concave, which implies $\Lambda = 0$, and where the optimization problem \eqref{eq:minimizing_scheme} is well-posed, the energy dissipation \eqref{eq:energy_stability} holds for every $\mathrm{dt} >0$. 
Therefore, the corresponding numerical scheme of the gradient flow~\eqref{eq:gradient_flow} is then \emph{unconditionally stable}.

\subsection{Deep Ritz method with energy splitting}
\label{subsec:DeepRitzSplit}

First introduced by Yu et al.~\cite{yu2018deep}, the \emph{Deep Ritz method} is a neural network-based approach for solving PDEs arising from variational problems. Rather than minimizing PDE residuals as in PINNs, it approximates the solution of a stationary variational problem by directly minimizing the corresponding energy functional. The method follows the same principle as the classical Ritz method, but replaces traditional basis functions by a neural network representation of the solution.

Namely, consider a PDE given by
\begin{equation}
  \nabla_{L^2} E(u) = 0  \,,
\end{equation}
where $E: V \to \mathbb{R}^+$  is an energy functional that takes non-negative values over a function space $V \subset L^2$.
The Deep Ritz method seeks a neural network $u_\vtheta$ such that
%
\begin{equation}
  u_{\vtheta} \approx \argmin_{v\in V} E(v).
\end{equation}
%


The classical Deep Ritz method is designed to approximate stationary solutions obtained as minimizers of an energy functional. In contrast, our objective is not to approximate a static minimizer but to learn the evolution operator of a dissipative dynamical system.
To achieve this, we exploit the variational formulation of the splitting scheme~\eqref{eq:minimizing_scheme}. Instead of minimizing the original free energy functional, we minimize the functional associated with a single time step of the energy-stable discretization. The resulting neural operator is therefore trained to approximate the time advancement operator defined by the splitting scheme rather than a stationary solution.
The free energy is decomposed into a convex and a concave part. At each time step, the concave contribution is linearized around the current state while the convex contribution is retained in its original form. This construction directly incorporates the variational structure of the splitting scheme into the training objective. As a result, the neural operator is trained to approximate an energy-stable time discretization of the underlying gradient flow.
Namely, we train a neural operator $\mathcal{NN}(\cdot; ; \vtheta)$ parameterized by a set of hyper-parameters $\vtheta$, such that for every $u_n \in V$, the output of the neural operator approximates the solution corresponding to the next time step:

\begin{equation}
  \begin{aligned}
	\mathcal{NN}(u_n; \vtheta) \simeq u_{n+1} 
= \argmin_{u \in V}&\bigg( 
\frac{1}{2 \mathrm{dt}} \| u - u_n \|_2^2 +
  E_1(u) 
                 \\ &\quad + E_2(u_n) + \langle \nabla_{L^2} E_2(u_n), u- u_n \rangle_{L^2}   \bigg).
  \end{aligned}
   \label{eq:deepRitzSplit}
\end{equation}

In Section~\ref{sec:AC} and Section~\ref{sec:dendritic_growth}, we will train the neural operators on the Allen-Cahn equation and the dendritic growth model, using the corresponding minimization scheme in its variational formulation 
directly as the loss function. In this way, the training procedure remains consistent with the dissipative variational structure of the underlying gradient flow.
We investigate its performance as function of the neural network architecture and compare it to conventional data-driven learning approaches.

Note that different approaches incorporating the energy dissipation property (\ref{eq:energy_dissipation}) can also be found in other works. For instance, in \cite{kutuk2024energy}, energy dissipation is enforced directly in the loss function as a penalty term, whereas in \cite{zhang2024energy}, a modified energy functional is introduced to ensure dissipation. However, the above works rely on training the neural operator with conventional data-driven approaches, whereas our training approach is unsupervised and requires only the initial conditions of the system.

\section{Application of the DeepRitzSplit method to the Allen-Cahn equation}
\label{sec:AC}

We consider here the Allen-Cahn equation \cite{allen1979microscopic}, a simple phase-field model that describes phase separation, defined on a domain $\Omega \subset \mathbb{R}^N$ for some $N\in\mathbb{N}$, by 
\begin{equation}
  \frac{\partial u}{\partial t} =  -\nabla_{L^2} E_\mathrm{AC}(u),\quad\text{on}\quad \Omega \times (0,T)
\label{eq:allen-cahn}
\end{equation}
which corresponds to a $L^2$-gradient flow of the Ginzburg-Landau free energy functional,
\begin{equation}
  E_\mathrm{AC}(u) = \int_\Omega \left( \frac{|\nabla u |^2}{2} + \frac{W(u)}{\epsilon^2} \right) d\vx \quad \forall u \in H^1(\Omega) \cap L^\infty(\Omega)
  \label{eq:CH-energy}
\end{equation}
$u: \Omega \times (0, T) \to \mathbb{R}$ is the phase-field solution of Eq.~\eqref{eq:allen-cahn}, where $\{u = -1\}$ and $\{u = 1 \}$ represent the two phases, respectively, $W(s) = (s^2 -1)^2/4$ is the double-well potential enforcing the phase separation, and $\epsilon > 0 $ characterizes the interface thickness. We also choose the initial condition $u(\cdot, 0) = u_0 \in H^1(\Omega) \cap L^\infty(\Omega)$ such that $\| u_0 \|_\infty \leq 1$.
A strong convergence result, proven by Evans et al.~\cite{evans1992phase}, states that the Allen-Cahn equation approximates the \emph{mean curvature flow} for a specific class of functions, which we will detail and use for benchmarks in Section~\ref{subsec:AC_benchmark}.

\subsection{DeepRitzSplit neural operator for the Allen-Cahn equation}
\label{subsec:AC_DRSscheme}

The first convex-concave splitting scheme for the Allen-Cahn equation was introduced by Eyre \cite{eyre1998unconditionally}, who decomposed the energy $\EAC$ as the sum of a convex part, $E_1$, and a concave part, $E_2$, defined by 
\begin{equation}
	\label{eq:AC-splitting-E1-E2}
	E_1(u) = \int_\Omega \left (\frac{|\nabla u |^2}{2} + \beta \frac{u^2}{2\epsilon^2} \right )d\vx
	\quad \text{and}\quad
	E_2(u) = \int_\Omega \left (\frac{W(u)}{\epsilon^2} - \beta \frac{u^2}{2\epsilon^2} \right )d\vx,
\end{equation}
for every $u \in H^1(\Omega) \cap L^\infty(\Omega)$ and for some $\beta >2$. Indeed, by the maximum principle, $\|u(\cdot, t) \|_\infty \leq 1$ for all $t\in (0, T)$, and hence $| W''(u)| \leq 2$. Consequently, $E_2''(u) \leq (2 - \beta) |\Omega|/ \epsilon^2 < 0$, which shows that $E_2$ is concave for $\beta >2$.

The semi-implicit numerical scheme is given by 
\begin{equation}
  \begin{aligned}
\frac{u_{n+1} - u_n}{\mathrm{dt}} &= -\nabla_{L^2} E_1(u_{n+1}) - \nabla_{L^2} E_2(u_n)
\\
& = \left (\Delta u_{n+1} - \frac{\beta}{\epsilon^2} u_{n+1}\right ) + \left(\frac{\beta}{\epsilon^2} u_n - \frac{W'(u_n)}{\epsilon^2}\right) \,,
  \end{aligned}
\label{eq:splitting_scheme_AC}
\end{equation}
which gives 
\begin{equation}
  u_{n+1} = \left( \mathrm{Id}_V - \mathrm{dt}\left(\Delta - \frac{\beta}{\epsilon^2} \right)\right)^{-1} \left( u_n - \frac{\mathrm{dt}}{\epsilon^2}(W'(u_n) - \beta u_n)\right) \,,
\label{eq:splitting_scheme_AC-sol}
\end{equation}
where $\mathrm{Id}_V$ is the identity operator on $V$ and $\Delta$ the Laplacian operator. The concavity of $E_2$, 
ensures that the splitting scheme~\eqref{eq:splitting_scheme_AC-sol} is unconditionally stable.

To establish the DeepRitzSplit method we construct a neural operator that represents a semi-group $\mathcal{NN}_{\theta, \mathrm{dt}}$ (or simply $\mathcal{NN}_\theta$ for the reader's convenience) of the isotropic Allen-Cahn equation, such that for every $u_n \in H^1 \cap L^\infty(\Omega)$ and a given time step $\mathrm{dt} > 0$,
\begin{equation}
   \mathcal{NN}(u_n ; \vtheta) \approx u_{n+1}
\end{equation}
where $u_{n+1}$ is defined by Eq.~\eqref{eq:splitting_scheme_AC-sol}.


The neural operator is trained using the variational loss associated with the splitting defined by $E_1$ and $E_2$ in Eq.~\eqref{eq:AC-splitting-E1-E2}:
\begin{equation}
	\mathcal{NN}(u_{n}; \vtheta) \approx \argmin_{u \in H^1(\Omega) \cap L^\infty(\Omega)} \mathcal{L}(u, u_n), \, \quad \forall u_n \in H^1(\Omega) \cap L^\infty(\Omega) \,,
	\label{eq:AC-NN-Ritz}
\end{equation}
where
\begin{equation*}
  \mathcal{L}(u,v) =   \frac{\| u - v \|^2_2}{2\mathrm{dt}}
                     + E_1(u) + E_2(v) 
                     + \langle \nabla_{L^2} E_2(v), u - v \rangle_{L^2(\Omega)}  \,.
\end{equation*} 
\deleted[id=MZ]{Since $\beta>2$, $E_2$ is concave and the underlying splitting scheme is unconditionally energy stable.} The trained neural operator approximates one step of \replaced[id=MZ]{the unconditionally energy-stable scheme~\eqref{eq:splitting_scheme_AC-sol}.}{this energy-stable scheme through the associated variational formulation.}
%
%

In practice, we train the neural operator on a finite set of initial conditions, which leads to minimizing the following empirical loss:
\begin{equation}
  \frac{1}{|K|}\sum_{u \in K}\mathcal{L}( \mathcal{NN}(u; \vtheta), u) \,,
\label{eq:AC-NN-Ritz-loss}
\end{equation}
where $K$ is the training configuration set and is a finite subset of $H^1(\Omega) \cap L^2(\Omega)$. 

\subsection{Reaction-Diffusion Neural Operator architecture} \label{subsec:RDNO}



For the experiments below, $\mathcal{NN}_\theta$ is implemented using the Reaction-Diffusion Neural Operator (RDNO) architecture shown in Fig.~\ref{fig:RDNO}. This architecture is inspired by the structure of the splitting scheme~\eqref{eq:splitting_scheme} and mirrors its decomposition into reaction and diffusion contributions. The reaction block represents the local nonlinear transformation, whereas the diffusion block represents the regularizing action of the implicit operator.
\deleted[id=MZ]{
The RDNO is composed of two main blocks, as illustrated in Fig.~\ref{fig:RDNO}.
Both blocks are learned in the Allen--Cahn experiments, consistently with $u_{n+1}=K_\beta * \rho_\beta(u_n)$.
}

%
%
Namely, we can write 
\begin{equation}
  \mathcal{NN}(\cdot\; ; \vtheta) = D_\vtheta \circ R_\vtheta
	\label{eq:RDNO}
\end{equation}
where
\begin{itemize}
\item 
   $R_\vtheta$, called the \emph{reaction block}, consists of a stack of convolutional layers 
   \[
   R_\theta = \pi_\theta \circ \left(\sigma \circ C^L_\theta \circ \ldots \circ \sigma \circ C^1_\theta\right) \circ T_\theta (u)
   \]
   with $\sigma$ the activation function, $T_\theta$ the lifting layer from discretized function space $V_N$ of $V$ to the higher-dimensional latent space $(V_N)^w$ with $w\in\mathbb{N}$ the lifting dimension, $\pi_\theta$  the projection layer from $(V_N)^w$ to $V_N$ and $C^i_\theta(\cdot) = A^i_\theta * \cdot + B^i_\theta$ the convolutional layers. 
\item  
   $D_\vtheta$, called the \emph{diffusion block}, employs an architecture that depends on the convex part of the energy splitting.
\end{itemize}

\begin{figure}[H]
  \centering 
  \scalebox{1.25}{
  \begin{tikzpicture}[circle dotted/.style={dash pattern=on .2mm off 3mm,
                                         line cap=round}]
    \draw[->] (-1.5, 0) node[left=5pt]{$u_n$} -- ++(0.5, 0) ;
    \draw[dashed] (-1, -1) rectangle ++(5, 2);
    \draw[] (-0.9, -0.9) rectangle ++(0.4, 1.8);
    \draw (-0.7, 0) node[rotate=90]{Lifting};

    \draw (-0.5, 0) -- ++(0.1, 0);
    \draw[] (-0.4, -0.9) rectangle ++(0.6, 1.8);
    \draw (-0.1, 0) node[rotate=90]{Conv2d};
    
    \draw (0.2, 0) -- ++(0.1, 0);
    \draw[] (0.3, -0.9) rectangle ++(0.6, 1.8);
    \draw (0.6, 0) node[rotate=90]{Activation};

    \draw (0.9, 0) -- ++(0.1, 0);
    \draw[line width = 0.5mm,circle dotted] (1.1, 0) -- (1.9,0);

    \draw (2, 0) -- ++(0.1, 0);
    \draw[] (2.1, -0.9) rectangle ++(0.6, 1.8);
    \draw (2.4, 0) node[rotate=90]{Conv2d};

    \draw (2.7, 0) -- ++(0.1, 0);
    \draw[] (2.8, -0.9) rectangle ++(0.6, 1.8);
    \draw (3.1, 0) node[rotate=90]{Activation};
    
    \draw (3.4, 0) -- ++(0.1, 0);
    \draw[] (3.5, -0.9) rectangle ++(0.4, 1.8);
    \draw (3.7, 0) node[rotate=90]{Projection};

    \draw (2, 1) node[above] {Reaction};
    \draw (4, 0) -- ++(0.5, 0);
    \draw[dashed] (4.5, -1) rectangle ++(1, 2);
    \draw (5, 1) node[above] {Diffusion};
    \draw[] (4.7, -0.9) rectangle ++(0.6, 1.8);
    \draw (5, 0) node[rotate=90]{Conv2d};

    \draw[->] (5.5, 0) -- ++(0.5, 0) node[midway, right=5pt]{$\mathcal{NN}(u_n;\vtheta)$};
\end{tikzpicture}
}
  \caption{An illustration of the Reaction-Diffusion Neural Operator}
  \label{fig:RDNO}
\end{figure}

This neural operator architecture is inspired by the work of Bretin et al.~\cite{bretin2022learning}, who used similar architectures for the Allen-Cahn equation. However, they used a supervised data-driven approach to train neural operators using samples generated by a classical convex-concave splitting scheme.
It is also worth noting that in the approach of~\cite{bretin2022learning}, the order of the Reaction and Diffusion blocks is inverted (DRNO) and that a multi-layer perceptron (MLP) is used for the reaction block, instead of convolutional layers. The RDNO architecture proposed here follows the structure of the splitting scheme and therefore establishes a direct correspondence between the numerical discretization and the neural operator design. In addition, using convolutional layers in the reaction block leads to a significantly smaller number of trainable parameters than a comparable MLP-based implementation, resulting in faster training.

\subsection{Numerical experiments}
\label{subsec:AC_benchmark}


The numerical experiments are organized in two stages. We first compare the unsupervised DeepRitzSplit and supervised data-driven training strategies using the same RDNO architecture. We then use DeepRitzSplit with three different neural-operator architectures: the bespoke RDNO, UNet, and FNO. 

\paragraph{Reference solver}
In both comparisons, the reference solutions are computed using a Fourier spectral spatial discretization combined with the splitting scheme in Eq.~\eqref{eq:splitting_scheme_AC-sol}.
We
consider the periodic domain $\Omega=[0,1]_{\mathrm{per}}^2$. 
Under periodic boundary conditions, the numerical scheme~\eqref{eq:splitting_scheme_AC-sol} can be written as
\begin{equation}
  u_{n+1}=K_\beta * \rho_\beta(u_n),
  \label{eq:AC-sol}
\end{equation}
where $K_\beta$ is the periodic convolution kernel defined by its Fourier
coefficients
\begin{equation}
  \widehat{K_\beta}(\xi) = \frac{1}{ 1+\mathrm{dt}\left(4\pi^2|\xi|^2+\beta\epsilon^{-2}\right) },
  \qquad \xi\in\mathbb{Z}^2,
  \label{eq:AC-kernel}
\end{equation}
where, for a periodic function $f$ on $[0,1]^2$,
\begin{equation}
  \widehat f(\xi) = \int_{[0,1]^2} f(x)e^{-2\pi\mathrm{i}\xi\cdot x}\,d x.
\end{equation}
Moreover,
\begin{equation}
  \rho_\beta(s) = s-\frac{\mathrm{dt}}{\epsilon^2} \left(W'(s)-\beta s\right).
\end{equation}


\paragraph{Training configurations}
%
We consider phase-field initial conditions whose zero level sets are \emph{connected} perturbed circles.
The behavior of the Allen-Cahn equation has been shown~\cite{evans1992phase} to be such that the $0$-level set of the corresponding solution converges to a curve evolving according to \emph{curve-shortening flow}, or \emph{mean curvature flow}, which decreases the perimeter of the interface during evolution.

Specifically, we consider initial conditions, $u_0$, parameterized by
\begin{equation}
u_0(x) = \tanh\left(
  \frac{r(\theta) - \| x - (0.5, 0.5)\|_{\text{per}}}{\sqrt{2}\epsilon}
\right), \quad \forall (x, y) \in [0, 1]_{\text{per}}^2 \,,
\label{eq:perturbed_disks}  
\end{equation}
where $\theta$ is the polar angle of the vector $x -(0.5,0.5)$, $\| \cdot \|_{\text{per}}$ is the canonical periodic distance on $[0, 1]^2_{\text{per}}$, and
\begin{equation}
 r(\theta) = r + r_p \sum_{k=1}^M \left[ \frac{a_k}{2 \sqrt{\sum_{j=1}^M a_j^2}} \cos (k\theta) +\frac{b_k}{2 \sqrt{\sum_{j=1}^M b_j^2}}\sin(k\theta) \right]\,,
\end{equation}
with $M$ the maximal mode of perturbations and $a_k, b_k$ sampled uniformly from $[-1, 1]$. 


%
\added[id=CK]{
  For the training of the neural operators, we consider a very small set of training configurations. These consist of 180 randomly generated samples of perturbed disks with a fixed parameter $M=5$. For each perturbed disk, we select a base radius $r$ from a uniform distribution $r\in[0.12,0.3]$, and for each chosen $r$, we select a perturbation radius $r_p$ from a uniform distribution $r_p \in[0.1r,0.5r]$.
}
%
For model validation, we used 60 samples with $M$, $r$, and $r_p$ within the parameter range of the training set (labeled as \emph{in-dist.}). We further used three types of \emph{out-of-distribution (OOD) configurations}:  
%

\begin{itemize}
  \item 
  Perturbed disks generated from Eq.~\eqref{eq:perturbed_disks} with $M=8$, with $r$ sampled uniformly from $[0.06,0.375]$ and, conditionally on $r$, $r_p$ sampled uniformly from $[0.09r,0.55r]$.
  \item 
  Configurations consisting of the union of three disks. Their radii $r_i$ are sampled independently and uniformly from $[0.12,0.3]$, and their centers $c_i$ are sampled independently and uniformly from $[0,1]^2_{\mathrm{per}}$. Overlapping disks are allowed, and the corresponding phase field is defined by
  \[
    u_0(x) = \max_{i=1,2,3}
               \left\{
                 \tanh\left( \frac{r_i-\lVert x-c_i\rVert_{\text{per}}}{\sqrt{2}\epsilon} \right)
               \right\} \,.
  \]
  \item 
  Random fields obtained by sampling the value of $u_0$ independently at each grid point from the uniform distribution on $[-1,1]$.
\end{itemize}

%
%
Each out-of-distribution test set contains 500 independently generated configurations.

\begin{figure}[H]
  \begin{center}
 \includegraphics[width=\linewidth]{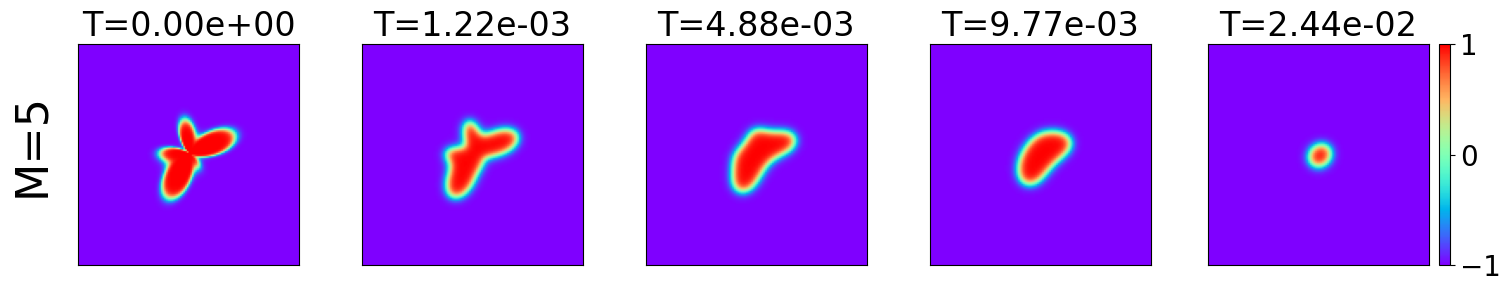} 
 \\
 \includegraphics[width=\linewidth]{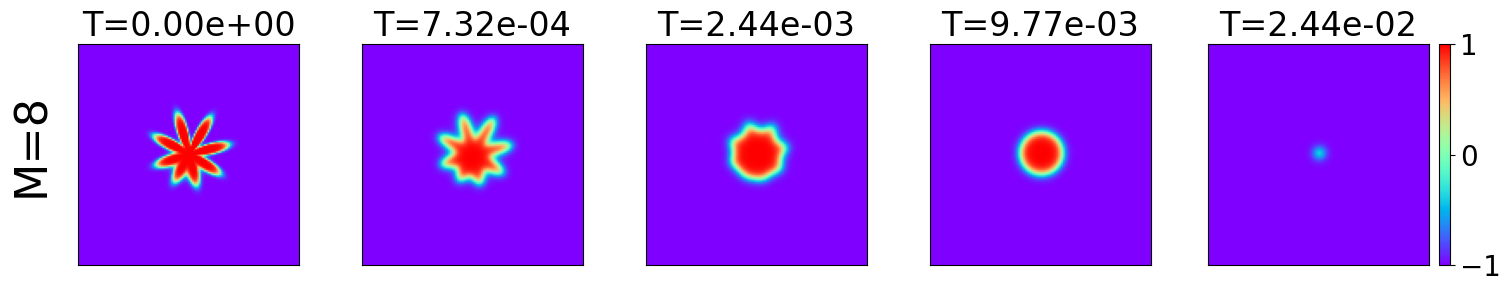} 
\end{center}
\caption{Perturbed disks evolving under mean curvature flow with $M=5$ (top) and $M=8$ (bottom). The numerical simulations are done with Fourier spectral method employing the splitting scheme defined in~\eqref{eq:splitting_scheme_AC-sol}.}
 \label{fig:perturbed_disks}
\end{figure}

The parameters for the benchmark are given in Table~\ref{tab:ac_params},
and
examples of evolution of perturbed circles under mean curvature flow are shown in Fig.~\ref{fig:perturbed_disks}.

\begin{table}[H]
  \caption{Parameters for numerical experiments of the Allen-Cahn equation}
 \label{tab:ac_params}
\centering
\begin{tabular}{cccc}
\hline
  \multicolumn{1}{l}{\bf $\beta$ }  &\multicolumn{1}{c}{\bf time step $\mathbf{\mathrm{dt}}$} &\multicolumn{1}{c}{\bf grid size}&\multicolumn{1}{l}{\bf Interface thickness $\epsilon$}
\\[1mm] 
$2.001$ & $\num{2.44e-4}$ & $128\times128$ & $1/64$
\\
\hline \\
\end{tabular}
\end{table}

All neural operators are trained to approximate the one-step mapping $u_n\mapsto u_{n+1}$ using only the 180 initial configurations $u_0$. In the unsupervised setting, the loss in Eq.~\eqref{eq:AC-NN-Ritz-loss} is evaluated directly on these configurations, without requiring reference solutions. In the supervised setting, each input $u_0$ is paired with a target $u_1$ generated \replaced[id=MZ]{with the Reference solver.}{by applying the splitting scheme~\eqref{eq:splitting_scheme_AC-sol}.}

%

\paragraph{Compared methods and architectures}
The first objective of the numerical experiments is to investigate the performance of the DeepRitzSplit method, particularly its accuracy, the ability to generalize to out-of-distribution configurations, as well as training and evaluation times. We compare the DeepRitzSplit unsupervised method, employing the loss function defined in Section~\ref{subsec:DeepRitzSplit}, to a conventional data-driven supervised neural method. The data-driven method is trained with solutions given by the Reference solver employing the numerical splitting scheme~\eqref{eq:splitting_scheme_AC-sol} and the loss function is the difference between the neural-operator prediction and the reference solution. 
For both 
traning strategies, we use the same
\emph{prescribed} RDNO architecture, as schematized in Fig.~\ref{fig:RDNO}.

The second objective is to evaluate the performance of the DeepRitzSplit with different neural architectures. We compare:
\begin{itemize}
	\item 
	An RDNO architecture, with a reaction block composed of a double convolutional layer of a width of $10$ and the hyperbolic tangent function for activation and a diffusion block consisting of a convolutional layer with a kernel size of $17$.
   \item 
   {Due to its popularity in biomedical image analysis, we consider }a UNet with 4 downsampling and upsampling levels, incorporating some modifications of the classical UNet by Falk et al.~\cite{falk2019u}: simple linear convolution instead of max pooling for downscaling, transposed convolutional layers for upscaling, downscaling and upscaling performed using a $3 \times 3$ convolution with stride 2, and periodic padding to align with the periodic boundary conditions in our configuration.
   \item 
   {Since the equation is endowed with periodic boundary conditions, it is also natural to consider a }Fourier Neural Operator (FNO)~\cite{Li2020FourierNO}, composed of 4 {{spectral}} layers, each with 20 modes and a width of 16.
\end{itemize}

\paragraph{Implementation and training details}
All models were trained with the Adam optimizer~\cite{kingma2014adam} using a batch size of $32$ and a learning rate of $10^{-3}$ on an NVIDIA RTX 6000 Ada GPU.
Training was stopped when the relative fluctuation of the training loss over the preceding $10000$ epochs fell below $3\text{\textperthousand}$. Across all models, the training time was approximately 5 minutes, with only minor variations between runs.
%



\paragraph{Evaluation Metrics}
In order to assess the accuracy of the neural operator solutions we compare the error fields as the absolute difference to the splitting method reference results,
i.e.
$e_n(x) = \left| u_{\mathrm{NO},n}(x)-u_{\mathrm{ref},n}(x) \right|$.
Furthermore we compare the shapes of the evolving phase, a key feature of the solution, given by the $u=0$ level set of the fields. 
In the cases with perturbed disks and multi-disks we additionally use the perimeter of the $u=0$ level set as a global metric for comparison with the reference splitting scheme solution. The perimeter, $P^\epsilon$, is a fundamental metric for the mean curvature flow given by the Allen-Cahn equation (see Modica-Mortola~\cite{de1977gamma}) and is defined by
\begin{equation}
  P^\epsilon(u) := \int_\Omega \left (\epsilon \frac{|\nabla u |^2}{2} + \frac{W(u)}{\epsilon}\right ) d\vx, \quad \forall u \in H^1(\Omega)\cap L^\infty(\Omega)
  \label{eq:perimeter}
\end{equation}
We then compute the average of relative errors 
\[ \frac{\left \vert P^\epsilon(u_{\text{neural operator}}) - P^\epsilon(u_{\text{splitting}})\right \vert}{P^\epsilon(u_{\text{splitting}})}
\]
%
Using the splitting-scheme solution as the reference, we compute the relative error in $P^\epsilon$ over 50 time steps, corresponding to the final time $T=50$ with $\mathrm{dt}=\num{1.22e-2}$.
%

\paragraph{Results}


The principal results are summarized in Table~\ref{tab:ac_RDNO_UNET_rel_err}. The corresponding fields, interfaces, and error distributions are shown in Figs.~\ref{fig:ac-ritz-vs-data-random} and~\ref{fig:AC_OOD}. Additional results are provided in Appendix~\ref{app:NO_AC}.

\begin{table}[H]
  \centering
  \caption{Mean relative errors in $P^\epsilon$ for the different training methods and neural-operator architectures, using the splitting scheme as the reference.}
  \label{tab:ac_RDNO_UNET_rel_err}
  \small
  \setlength{\tabcolsep}{5pt}
  \begin{tabular}{llcccc}
    \hline
      & & \multicolumn{1}{c}{\bf In-distribution} & \multicolumn{3}{c}{\bf Out of distribution} \\
    \cline{3-3}\cline{4-6}
    \bf Method
      & \bf Architecture
      & \shortstack{Perturbed\\disks}
      & \shortstack{Perturbed\\disks}
      & \shortstack{Multiple \\disks}
      & \shortstack{Random   \\fields} \\
    \hline
    DeepRitzSplit  &  RDNO  &   0.16\%  &   0.48\%  &   0.37\%  &   1.02\%  \\
    Data-driven    &  RDNO  &   0.89\%  &   1.82\%  &   2.59\%  &  43.76\%  \\[1mm]
    DeepRitzSplit  &  UNet  &   0.45\%  &   0.28\%  &   1.07\%  &   2.49\%  \\
    Data-driven    &  UNet  &   1.02\%  &   1.29\%  &   2.67\%  &  67.05\%  \\[1mm]
    DeepRitzSplit  &  FNO   &  60.39\%  &  25.00\%  &  12.92\%  &  80.00\%  \\
    \hline
  \end{tabular}
\end{table}


\emph{Effect of the training strategy.}
We first compare DeepRitzSplit with supervised data-driven training using the same RDNO architecture. The two models are trained using the same 180 input configurations. As shown in Table~\ref{tab:ac_RDNO_UNET_rel_err}, the RDNO trained with DeepRitzSplit gives smaller errors than the \replaced[id=MZ]{data-driven}{supervised} RDNO for all four test sets. The difference is particularly large for the random fields, for which the \replaced[id=MZ]{data-driven}{supervised} RDNO does not reliably reproduce the reference evolution; see Fig.~\ref{fig:ac-ritz-vs-data-random}.
Although less accurate than for the perturbed-disk and multiple-disk configurations, the RDNO trained with DeepRitzSplit still reproduces the main features of the reference solution for this substantially different type of initial condition.

The results obtained with the UNet architecture, summarized in Table~\ref{tab:ac_RDNO_UNET_rel_err} and in Figs.~\ref{fig:AC_OOD}, \ref{fig:AC_IOD} and~\ref{fig:AC_multi}, show the same overall trend. For perturbed and multiple disks the error of the DeepRitzSplit UNet is somewhat smaller than with the data-driven method. For random fields the DeepRitzSplit UNet reproduces the main features while the data-driven method fails.

For the OOD example shown in Fig.~\ref{fig:AC_OOD}c-d, the errors of the \replaced[id=MZ]{data-driven}{supervised} models appear to be distributed relatively uniformly along the diffuse interface. Together with the large errors obtained for random fields, this suggests that the \replaced[id=MZ]{data-driven}{supervised} models have difficulty generalizing to interface configurations that differ substantially from those represented in the training set.

Using the same number of input configurations provides a controlled comparison in this limited-data setting. The results suggest that the \replaced[id=MZ]{data-driven}{supervised} models may require more training pairs to generalize to
configurations that differ substantially from the training distribution. They do not, however, establish that data-driven supervised training would necessarily remain less accurate with a larger training dataset.

\begin{figure}[H]
  \centering
  \includegraphics[width=\linewidth]{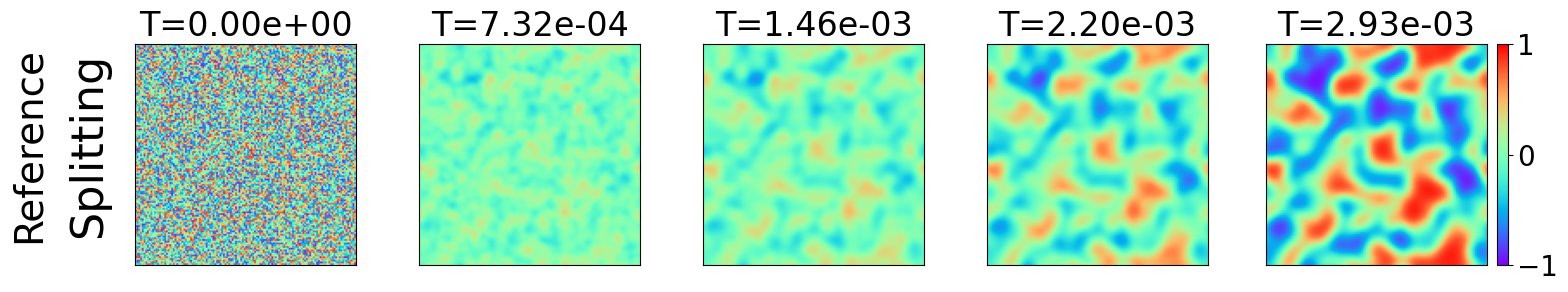}
  \includegraphics[width=\linewidth]{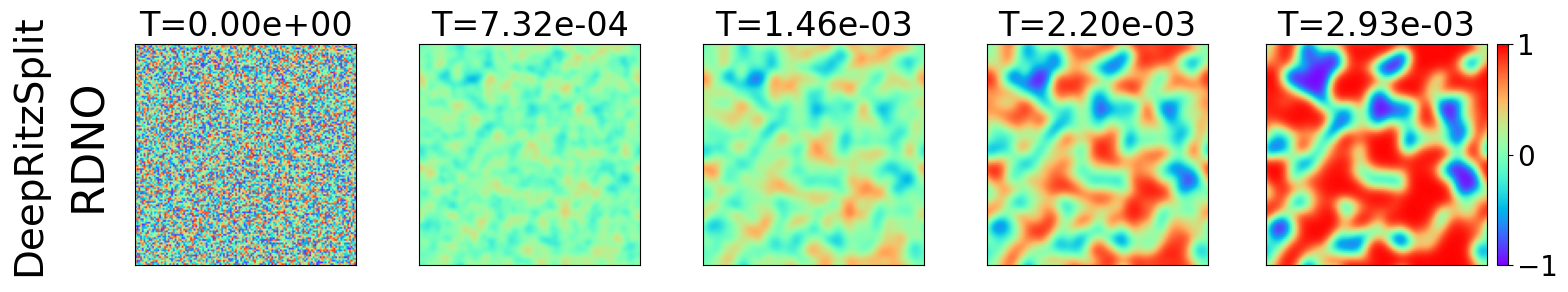}
  \includegraphics[width=\linewidth]{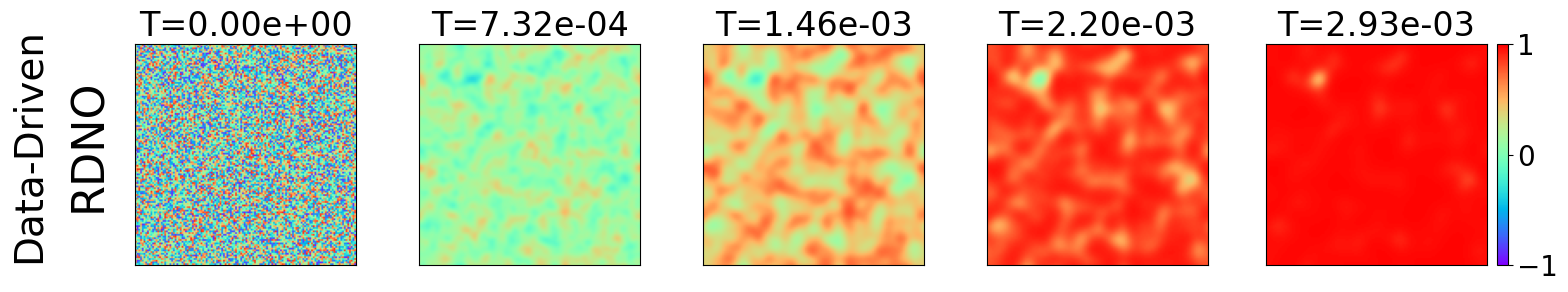}
  \caption{
   Evolution from a random-field initial condition. From top to bottom: reference solution computed with the splitting scheme, prediction of the RDNO trained with DeepRitzSplit, and prediction of the RDNO trained with supervised data-driven learning.
      }
  \label{fig:ac-ritz-vs-data-random}
\end{figure}


\emph{Effect of the neural-operator architecture.}
We next compare RDNO, UNet, and FNO within the DeepRitzSplit framework (Table~\ref{tab:ac_RDNO_UNET_rel_err}). RDNO gives the smallest errors for the in-distribution perturbed disks, the multiple-disk configurations, and the random fields. UNet gives the smallest error for the OOD perturbed disks. Thus, RDNO performs better for three of the four test families, but it does not systematically outperform UNet.

For the example shown in Fig.~\ref{fig:AC_OOD}, the $u=0$ level sets predicted by RDNO and UNet remain close to the reference interfaces over the displayed time interval. Their error fields are concentrated mainly around the diffuse interface, but their spatial distributions differ. With RDNO, the largest discrepancies tend to occur in regions of higher interface curvature. With UNet, the error is less clearly localized and appears larger in magnitude for this example.

Since $P^\epsilon(u)=\epsilon E_{\mathrm{AC}}(u)$ for the Allen-Cahn model, the evolution of $P^\epsilon$ also provides information about the energy evolution. For the example shown in Fig.~\ref{fig:AC_OOD}, the DeepRitzSplit predictions remain closer to the reference evolution than the supervised predictions. Among the DeepRitzSplit models, the UNet prediction decreases faster than the splitting-scheme reference, whereas the RDNO remains closer to the reference over the displayed time interval.

In the present DeepRitzSplit experiments, the FNO architecture produces substantially larger errors than with RDNO and UNet and does not reproduce the interface shape accurately. One possible explanation is the spectral truncation used in the FNO. The evolution of the thin diffuse interface may involve frequencies beyond the 20 modes retained in each spectral layer, which could limit the accuracy of the learned one-step operator. This interpretation remains speculative, since the results may also depend on the selected number of modes and other architectural and training parameters. Additional experiments would be needed to determine the origin of the observed performance difference.

\begin{figure}[H]
  \centering
  \includegraphics[width=.8\textwidth]{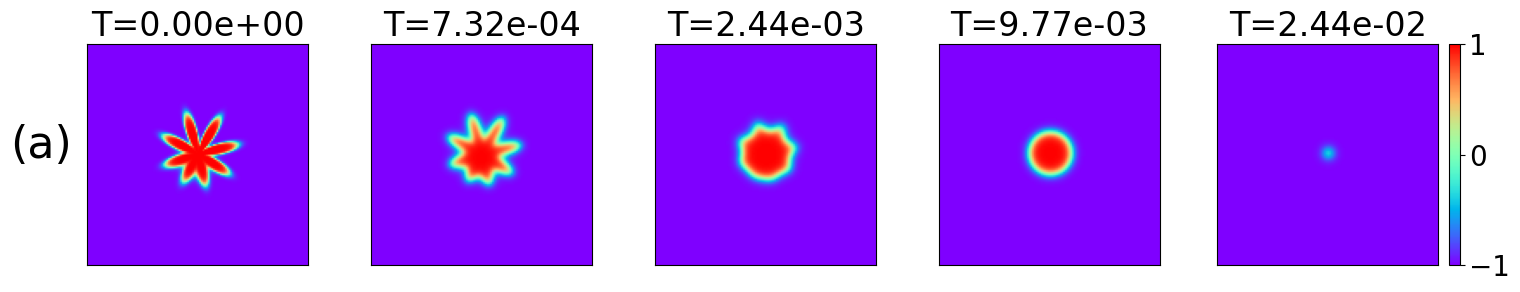}
  \includegraphics[width=.8\textwidth]{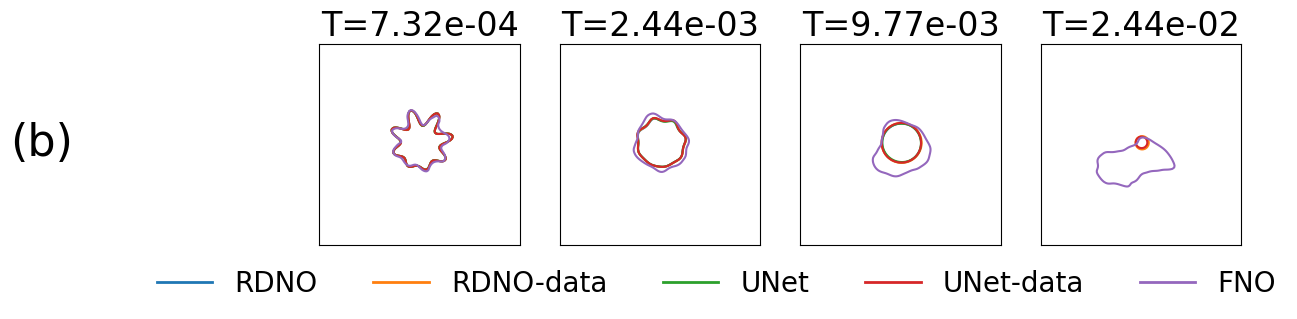}
  \includegraphics[width=.8\textwidth]{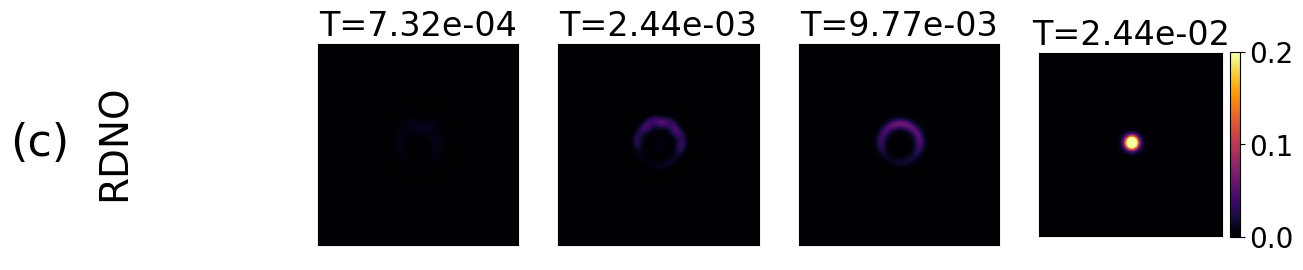}
  \includegraphics[width=.8\textwidth]{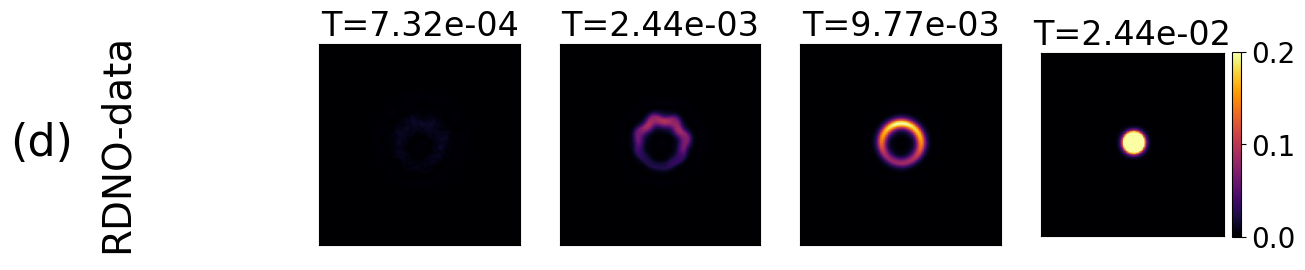}
  \includegraphics[width=.8\textwidth]{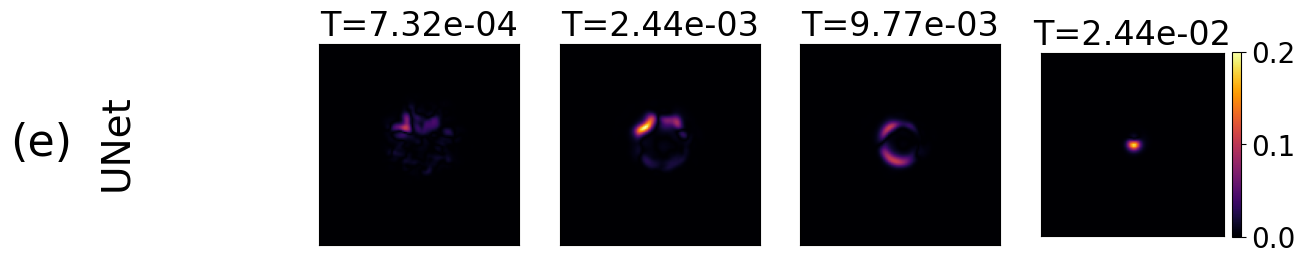}
  \includegraphics[width=.8\textwidth]{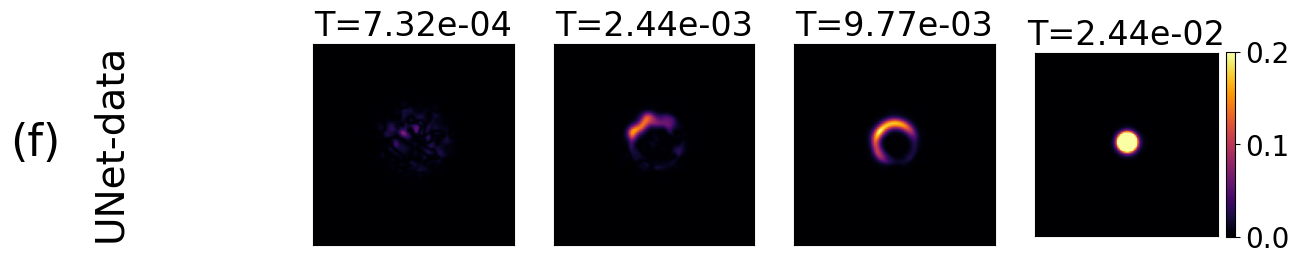}
  \includegraphics[width=.8\textwidth]{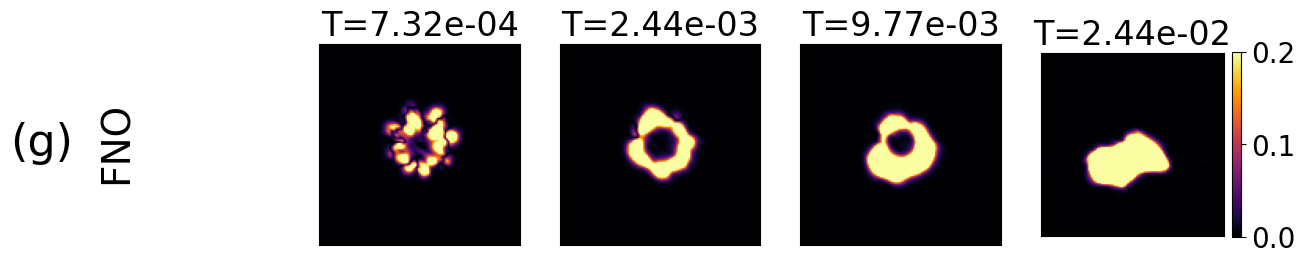}
  \includegraphics[width=.8\textwidth]{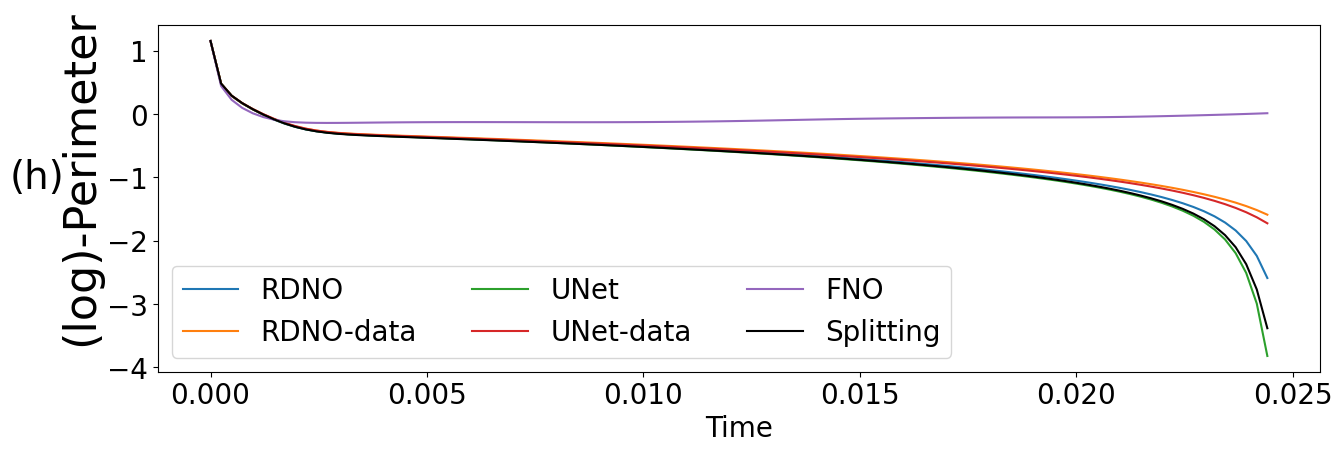}
  \caption{
      Predictions for an OOD perturbed-disk configuration.
      (a) Reference solution computed with the splitting scheme.
      (b) Comparison of the predicted interfaces through their
      $u=0$ level sets.
      (c--g) Absolute error fields for DeepRitzSplit with RDNO,
      supervised training with RDNO, DeepRitzSplit with UNet,
      supervised training with UNet, and DeepRitzSplit with FNO,
      respectively.
      (h) Evolution of the diffuse-interface perimeter $P^\epsilon$.
         }
  \label{fig:AC_OOD}
\end{figure}

\emph{Computational cost.}
The numbers of trainable parameters and the evaluation times are reported in Table~\ref{tab:ac_RDNO_UNET_eval_time}. Since evaluation time depends on the architecture rather than on the training strategy, the supervised and DeepRitzSplit versions of a given architecture have the same evaluation cost. For 50 time steps, the Fourier splitting method (Reference solver) is faster than all three neural operators in this Allen--Cahn test case. The comparison is different for the more complex dendritic-growth application considered in Section~\ref{sec:dendritic_growth}.

\begin{table}[H]
  \centering
  \caption{Number of trainable parameters and evaluation time over
  50 time steps for the neural-operator architectures and the reference
  splitting scheme.}
  \label{tab:ac_RDNO_UNET_eval_time}
  \begin{tabular}{lcc}
    \hline
    Method
      & Trainable parameters
      & Evaluation time (s) \\
    \hline
    RDNO      & 31522 & 0.181 \\
    UNet      & 37021 & 0.263 \\
    FNO       & 82174 & 0.320 \\
    Splitting & n/a   & 0.124 \\
    \hline
  \end{tabular}
\end{table}
\section{Neural operators for anisotropic dendritic growth}
\label{sec:dendritic_growth}

Compared to the studies on deep learning models for the isotropic Allen-Cahn equation, the use of neural networks to approximate phase-field models of anisotropic dendritic growth during solidification remains relatively unexplored. This is primarily due to the added mathematical complexity introduced by the anisotropy term.
In the scarce prior work~\cite{tseng2023deep, oommen2024rethinking} exclusively supervised data-driven approaches are employed. In contrast, we use the unsupervised physics-informed DeepRitzSplit neural operator to address the problem. Furthermore, we provide a comprehensive analysis of the crucial characteristics of predicted dendritic morphology, including the dendrite tip velocities and tip radii.

We consider the growth of dendritic crystal grains from a pure liquid using the well-established quantitative phase-field model of Karma \& Rappel~\cite{karma1998quantitative}. The variational formulation of the model is given by the system of equations
 \begin{equation}
\begin{cases}
\tau \phi_t &=  -\nabla_\phi E(\phi, U)
\\ 
U_t &= D \Delta U + K h'(\phi) \phi_t \,,
\end{cases}
\label{eq:dendritic_growth}
\end{equation}
where $\phi$ is the phase of the material and $U$ is the \emph{dimensionless} temperature, such that $U=0$ is the dimensionless melting temperature. $\tau>0$ is the relaxation time, 
$D > 0 $ is the diffusion coefficient and  $K>0$ is the latent heat parameter.
The free energy of the system is expressed as 
\begin{equation}
E(\phi, U) = \int_\Omega \left ( 
    \frac{1}{2} | a(\phi)^2| |\nabla \phi|^2 + \frac{\lambda_0}{2\epsilon K} U^2 + \frac{W(\phi)}{\epsilon^2}  +  \frac{\lambda_0}{\epsilon} h(\phi) U
\right )\, dx
\label{eq:free-energy}
\end{equation}
with 
$$ 
h(s) := s^5/5 - 2s^3/3 +s \,.
$$
$W=h'/4$ is the double-well function and $\epsilon>0$ is the interface thickness. The anisotropy of the interface energy is given by $a(\phi) = 1 + \sigma \cos(m \vec \Theta)$, where 
$\Theta = \arctan(n_2/n_1)$ 
is the in-plane azimuthal angle with the normal field $(n_1,n_2) = \vec n = \nabla \phi / |\nabla \phi |$ that specifies the orientation of the interface relative to the crystalline axis, which is in our case aligned with the horizontal axis. $m \in \mathbb{N}$ is the symmetry order of the crystallization, $\sigma>0$ is the anisotropy strength. Finally, we choose 
$\lambda_0  = \displaystyle\frac{D \tau}{0.6267 \epsilon}$,
in order to correctly reproduce negligible interface kinetics in the thin-interface limit~\cite{karma1998quantitative}.

\begin{figure}[ht]
	\centering
	\includegraphics[width=\linewidth]{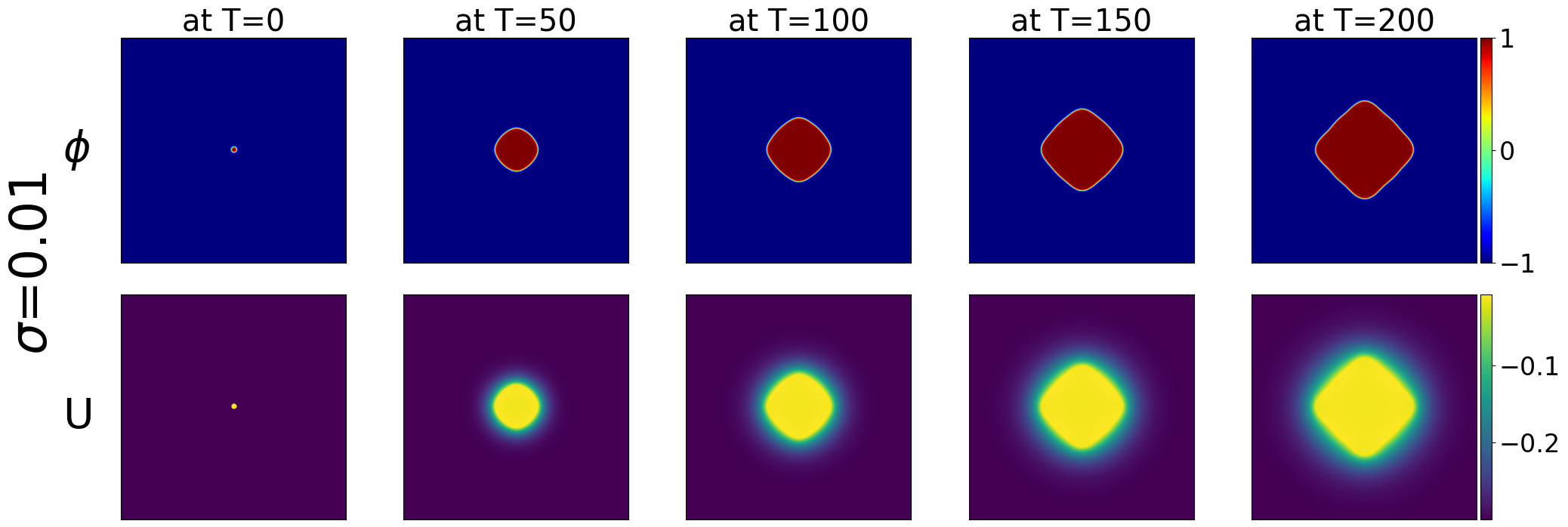}
	\\
	\includegraphics[width=\linewidth]{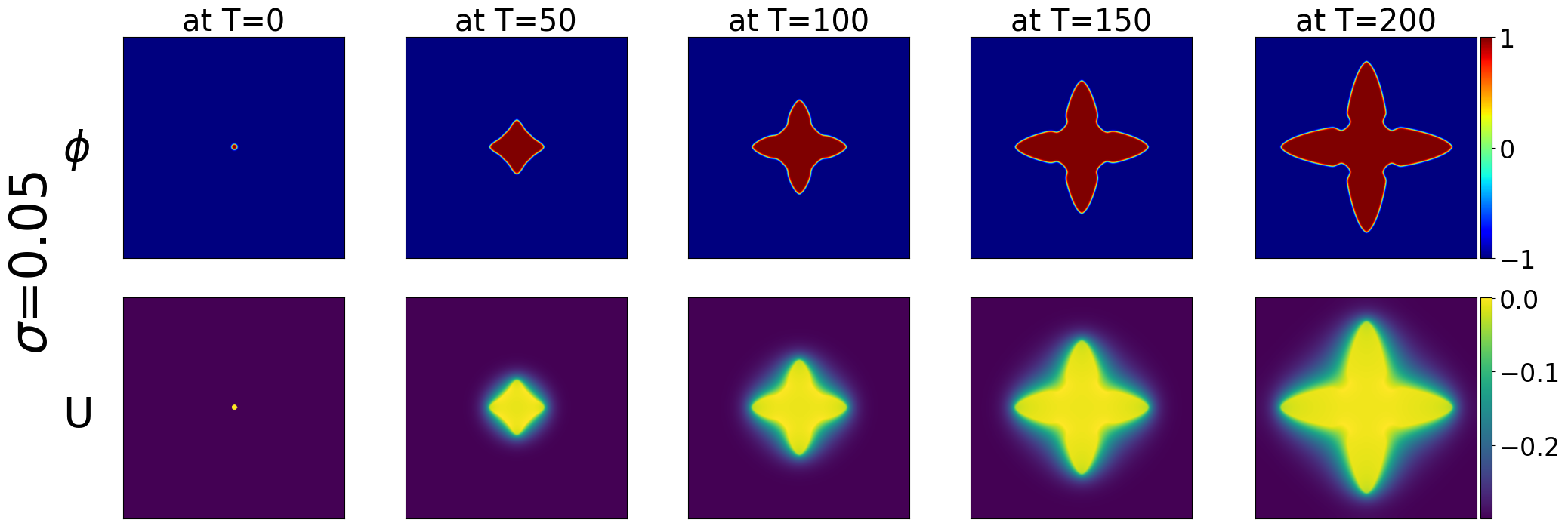}
	\caption{Dendritic growth simulations with fourfold symmetry ($m=4$) and anisotropy strengths $\sigma = 0.01$ (top) and $\sigma= 0.05$ (bottom) on a unit square. Simulation with the SAV method.}
	\label{fig:dendrite_with_diff_aniso}
\end{figure}

Two examples of dendrite growth for $\sigma=0.01$ and $\sigma=0.05$ are shown in Fig.~\ref{fig:dendrite_with_diff_aniso}, illustrating the growth evolution and the effect of the anisotropy strength, $\sigma$. 
\footnote{We restrict our analysis to cases where $\sigma < 1/15 \approx 0.067$. For $\sigma \geq 1/15$, the free energy $E$ defined in Eq.~\eqref{eq:free-energy} no longer remains convex, leading to the formation of \emph{faceted dendrites}. This phenomenon lies outside the scope of the present paper.}
The initial conditions for these and all following simulations are an initial solid crystal in form of a disk of radius $5\epsilon$ and a uniform initial temperature of $U_0=\kappa=-0.3$ in the liquid and $U_0=0$ in the solid. The boundary conditions are periodic.

\subsection{DeepRitzSplit formulation for dendritic growth}\label{subsec:dendritic_DRS}


In order to apply the DeepRitzSplit method (given in Section~\ref{subsec:DeepRitzSplit}) to the phase-field model, Eq.~\eqref{eq:dendritic_growth}, we split the free energy $\Ed$, defined by Eq.~\eqref{eq:free-energy}, into two parts:
$$ 
\Ed(\phi , U) = E_1(\phi, U ) + E_2(\phi, U) \,,
$$ 
where 
\begin{equation}
E_1(\phi, U) = \int_\Omega \left ( 
    \frac{1}{2} a(\phi)^2 |\nabla \phi|^2+ \beta \frac{\phi^2}{2 \epsilon^2} 
\right ) d\vx
\label{eq:E1-dendritic-growth}
\end{equation}
and 
\begin{equation}
E_2(\phi, U) 
=\int_\Omega \left (\frac{\lambda_0}{2\epsilon K} U^2 + \frac{W(\phi)}{\epsilon^2} + \frac{\lambda_0}{\epsilon}h(\phi) U  - \frac{\beta}{2\epsilon^2} \phi^2 \right) d\vx
\label{eq:E2-dendritic-growth}
\end{equation}
for some $\beta >0$.


Then, the corresponding solution scheme, employing the semi-implicit splitting scheme for the phase-field equation and a fully implicit scheme for the heat diffusion equation, is written as:
\begin{equation}
  \left \{
\begin{aligned}
 \tau \frac{\phi_{n+1}- \phi_n}{\mathrm{dt}} &=  -\nabla_\phi E_1(\phi_{n+1}, U_n) - \nabla_\phi E_2(\phi_n, U_n)
  \\ 
 \frac{U_{n+1} - U_n}{\mathrm{dt}} &= D \Delta U_{n+1} + K h'(\phi_{n+1}) \frac{\phi_{n+1} - \phi_n}{\mathrm{dt}}.
\end{aligned}
\right .
  \label{eq:aniso_scheme}
\end{equation}
with initial conditions $\phi_0, U_0 \in H^1(\Omega)\cap L^\infty(\Omega)$. For physical consistency, we also impose $\|\phi_0\|_\infty \leq 1$ and $\kappa \leq U_0 \leq 0$ with $\kappa$ the initial dimensionless temperature of the liquid.

\replaced[id=MZ]
{The first equation in~(\ref{eq:aniso_scheme}) is the semi-implicit splitting scheme for the phase-field equation, (see also Section~\ref{subsec:stable_scheme}), which we solve through the corresponding minimization problem:}
{The first phase-field equation in~\eqref{eq:aniso_scheme} is then solved with the semi-implicit splitting scheme~(\ref{eq:splitting_scheme}), which leads to the following minimization problem:}
\begin{equation}
  \begin{aligned}
\phi_{n+1} 
:= \argmin_{\phi \in V}&\bigg( 
\frac{\tau}{2 \mathrm{dt}} \| \phi - \phi_n \|_2^2 +
  E_1(\phi, U_n) 
                 \\ &\quad + E_2(\phi_n, U_n) + \langle \nabla_\phi E_2(\phi_n, U_n), \phi- \phi_n \rangle_{L^2}   \bigg).
  \end{aligned}
\label{eq:aniso_splitting}
\end{equation}
\replaced[id=MZ]
{Whereas the second equation in~\eqref{eq:aniso_scheme} is the implicit scheme for the temperature equation, which we solve in the discretized Fourier space, through the scheme}
{Whereas the second equation is solved in the discretized Fourier space with a implicit scheme, given by}
\begin{equation}
  U_{n+1}= \left ( \Id - \mathrm{dt}D \Delta \right)^{-1} \left( U_n + \mathrm{dt} K h'(\phi_{n+1}) \frac{\phi_{n+1} - \phi_n}{\mathrm{dt}}\right) \,.
\end{equation}
Note that this minimization problem has a solution for all $\phi_n, U_n$ since the functional in Eq.~\eqref{eq:aniso_splitting} is $C$-elliptic with $C \geq (1-\sigma)^2/4 > 0$. By the maximum principle, $\|\phi(t) \|_\infty \leq 1$ and $ \kappa \leq U(t) \leq 0$ for all $t\geq 0$, which gives
\begin{equation}
 D_\phi^2 E_2 (v, v) \leq \frac{1}{\epsilon}(2+ 4 \lambda_0 \epsilon \kappa ( \sqrt{1/3}^3 - \sqrt{1/3}) - \beta) \| v \|_2 \,\quad \forall v \in L^2(\Omega) \,,
  \label{eq:concavity_E2}
\end{equation}
where $D^2_\phi$ is the second-order Gateaux differential at $\phi$ on $L^2(\Omega)$.
Therefore, with a well-chosen $\beta >0$, depending only on $\kappa$, $\lambda_0$, and $\epsilon$, $E_2$ becomes concave, so that the numerical scheme~\eqref{eq:aniso_scheme} is unconditionally stable and thus the energy dissipation property~\eqref{eq:energy_dissipation} is retrieved. 

The presence of the anisotropic term makes the implicit part of Eq.~\eqref{eq:aniso_scheme} substantially more difficult to solve than in the isotropic Allen--Cahn case.
This motivates the use of a neural operator to solve the underlying optimization problem \eqref{eq:aniso_splitting}. More precisely, we aim to construct a neural operator $\mathcal{NN}(\cdot\;;\vtheta)$ such that 
 \begin{equation}
   \mathcal{NN}(\phi_n, U_n;\vtheta) \approx\phi_{n+1} \,,
  \label{eq:PF-Ritz}
 \end{equation}
where $\phi_{n+1}$ is a minimizer of Eq.~\eqref{eq:aniso_splitting} and thus the solution of the scheme given by Eq.~\eqref{eq:aniso_scheme}.

As in the case of the Allen-Cahn equation, this leads to minimizing the following empirical loss: 
\begin{equation}
  \frac{1}{|K|} \sum_{(\phi_n, U_n)\in K} \mathcal{L}( \mathcal{NN}(\phi_n, U_n; \vtheta), \phi_n, U_n) \,,
  \label{eq:PF-Ritz-loss}
\end{equation}
where 
\[
  \begin{aligned}
  \mathcal{L}(\phi, \psi, U) :=
  \bigg( 
        \frac{\tau}{2 \mathrm{dt}} \| \phi - \psi \|_2^2 + E_1(\phi, U) 
         + E_2(\psi, U) + \langle \nabla_\phi E_2(\psi, U), \phi- \psi \rangle_{L^2}   
  \bigg)
  \end{aligned}
\]
and the training configuration set $K$ is a finite subset of $(H^1(\Omega)\cap L^\infty(\Omega))^2$.

For benchmarking the neural operator for dendritic growth we use the Dendritic Growth benchmark problem of the \href{https://pages.nist.gov/pfhub/benchmarks/benchmark3.ipynb/}{The Phase Field Community Hub}. All parameters are given in Table~\ref{tab:pf_params}. Note that the scaling of the spatial coordinates is different from the published formulation -- instead of the interface thickness, $\epsilon$, we use the length of the square domain as the characteristic length to scale the dimensionless coordinates. According to Eq.~\eqref{eq:concavity_E2}, the choice of $\beta=14.4$ ensures the concavity of $E_2$ since $ D_\phi^2 E_2(v, v) < \frac{1}{\epsilon}( 2+ (-19.14) \cdot (-0.385) - \beta ) \| v\|^2_2 \approx -\frac{1}{\epsilon} 5.03 \| v\|_2^2< 0$ for every $v \in L^2(\Omega) \; \backslash\; \{0 \}$.

  \begin{table}[H]
  \caption{Parameters for numerical experiments of the phase-field model of dendritic growth. Note that all parameters are dimensionless.}
 \label{tab:pf_params}
\begin{center}
\begin{tabular}{lcr}
	\hline
	Anisotropy strength   &   $\sigma$    &     $0.01$ or $0.05$ \\
	Symmetry order        &      $m$      &                  $4$ \\
	Interface thickness   &  $\epsilon$   &              $1/400$ \\
	Relaxation time       &    $\tau$     &   $1.6 \cdot 10^{5}$ \\
	Coupling parameter    &  $\lambda_0$  &                $6382.64$ \\
	Diffusion coefficient &      $D$      & $6.25 \cdot 10^{-5}$ \\
	Latent heat parameter &      $K$      &                $0.5$ \\
	Initial temperature   &   $\kappa$    &              $-0.3$ \\
	Splitting parameter   &    $\beta$    &               $14.4$ \\
	Time step             & $\mathrm{dt}$ &               $0.04$ \\
	Grid spacing          &      $h$      &              $1/400$ \\ \hline
\end{tabular}
\end{center}
\end{table}

\subsection{Prescribed-Reaction RDNO architecture}\label{subsec:PR-RDNO}


\added[id=MZ]{For the dendritic growth model we introduce a particular RDNO architecture, where the reaction operator is not approximated by a neural network but is rather calculated explicitly. The splitting scheme (\ref{eq:splitting_scheme}) can be written as}
%
\begin{equation}
  \mathcal{P} (u_{n+1}) =  \mathcal{R}(u_n) \,,
  \label{eq:split_op1}
\end{equation}
where 
\begin{equation}
  \mathcal{P}(u): = u + \mathrm{dt} \nabla E_1(u) \quad\text{and}\quad \mathcal{R}(u) = u - \mathrm{dt} \nabla E_2(u).
\label{eq:split_op2}
\end{equation}

If there exists a neural operator $\mathcal{D}_\vtheta$, such that 
\begin{equation}
  \mathcal{D}_\vtheta \circ \mathcal{P} \approx \mathrm{Id}_V,
  \label{eq:inverse_diffusion}
\end{equation}
the splitting scheme of Eq.~\eqref{eq:split_op1} develops into
\begin{equation}
  u_{n+1} \approx \mathcal{D}_\vtheta \circ \mathcal{R}(u_n).
  \label{eq:prescribed_react}
\end{equation}
The structure of the scheme leads us to consider the
\emph{Prescribed-Reaction RDNO} (PR-RDNO),
shown in Fig.~\ref{fig:RDNO_prescribed}, where $\mathcal{D}_\vtheta$ is the neural diffusion operator and $\mathcal{R}$ is the explicitly calculated (prescribed) reaction operator.

We use a UNet as diffusion block in the PR-RDNO architecture. Recall that, first introduced in \cite{ronneberger2015u}, UNet is a convolutional neural network architecture originally designed for biomedical image segmentation. It consists of a contracting path (encoder) that captures context through downsampling, and a symmetric expanding path (decoder) that enables precise localization through upsampling. Crucially, it includes skip connections between corresponding layers in the encoder and decoder to retain spatial information lost during downsampling. UNet became popular due to its ability to train on relatively few images while producing high-accuracy segmentations.

%
%

\begin{figure}[H]
  \centering 
  \scalebox{1.25}{
  \begin{tikzpicture}[circle dotted/.style={dash pattern=on .2mm off 3mm, line cap=round}]
    \draw[->] (-1.5, 0) node[left]{$u_n$} -- ++(.5, 0) ;
    \draw[dashed] (-1, -1) rectangle ++(4.2, 2);

    \draw (1.1, 1) node[above] {{\color{red}{Reaction}}};
    \draw (1.1, 0) node {{\color{red}{$\mathcal{R}(u_n)=u_n - \mathrm{dt} \nabla E_2(u_n)$}}};

    \coordinate (Blk2) at (3.2, 0);
    \draw (Blk2) -- ++(0.5, 0);
    \draw[dashed] (Blk2)++(0.5, -1) rectangle ++(1, 2);
    \draw (Blk2) ++(1, 1) node[above] {{\color{blue}{Diffusion}}};
    \draw[] (Blk2)++ (0.7, -0.9) rectangle ++(0.6, 1.8);
    \draw (Blk2)++(1, 0) node[rotate=90]{{\color{blue}UNet}};
    \draw[->] (Blk2) ++(1.5, 0) -- ++(0.5, 0); 
    \draw (Blk2) ++ (2,0) node[right=3pt] {$\mathcal{NN}(u_n;\vtheta)$};
\end{tikzpicture}
}
\caption{Prescribed-Reaction RDNO (PR-RDNO) with a UNet diffusion operator.}
  \label{fig:RDNO_prescribed}
\end{figure}


\subsection{Numerical experiments}

\paragraph{Reference solver}
The reference for comparison of all neural methods is the Scalar Auxiliary Variable (SAV)-based method for the Karma-Rappel phase-field model of anisotropic dendritic growth, which was recently developed by Guo et al.~\cite{guo2024efficient}. We use Fourier spectral spatial discretization for the SAV method. For the reader's convenience, technical details on the SAV scheme are described in Appendix~\ref{app:sav}.

\paragraph{Training configurations}
We consider dendrites growing in a square unit domain from one, two or three initial small disk-shaped nuclei of radius of $5 \epsilon$ into a liquid that is initially uniformly undercooled at $U=\kappa$. The initial condition is thus
\begin{equation}
	\label{eq:dendrite_initialcondition}
	\begin{aligned}
		\phi_0(x, y) &= \max_{i = 1, \ldots, n} \left \{ \tanh \left ( 
        \frac{5\eps - \|(x, y) - c_i\|_{\text{per}}}{\sqrt{2}\epsilon}
		\right )
		\right \}
		\quad\text{and}
		\\
		U_0 &=
		\begin{cases}
			0 &\text{ on } \{\phi_0 > 0 \}
			\\
			\kappa &\text{ otherwise }.
		\end{cases}
	\end{aligned}
\end{equation}
with $n \in \{1, 2, 3\}$ and $c_i \in [0, 1]^2_{\text{per}}$.
For the training of the neural operators we use 6,000 sample configurations for each anisotropic strength, $\sigma$. These sample configurations result from 12 different initial conditions: 4 randomly generated arrangements for each, \textbf{one}, \textbf{two}, and \textbf{three} grains in the domain, respectively, and 500 sample configurations from the evolution corresponding to each initial condition.

\added[id=CK]{
  In contrast to the training procedure used for the Allen-Cahn equation in Section~\ref{subsec:AC_benchmark}, the neural operator for dendritic growth is trained using all 6,000 of $(\phi_n, U_n)$ configurations and not only the initial condition $(\phi_0, U_0)$.  
  \replaced[id=MZ]
  {For the data-driven supervised training we use solutions given by the SAV method as data and the loss function is the difference between the neural-operator prediction and the SAV solution. The training configurations are the same as for the unsupervised training. For both methods the identical PR-RDNO architecture is used, as schematized in Fig.~\ref{fig:RDNO_prescribed}.}
  {For the data-driven supervised training, the network minimizes the discrepancy between its prediction and the next time-step solution $\phi_{n+1}$, which is generated by applying the SAV scheme~\eqref{eq:SAV_scheme} to the current state $(\phi_n, U_n)$.}
}

\paragraph{Compared methods and architectures}


The first objective of the numerical experiments is to investigate the performance of the DeepRitzSplit method, particularly its accuracy, the ability to generalize to out-of-distribution configurations, as well as training and evaluation times. We compare the DeepRitzSplit unsupervised method, employing the loss function defined in Eq.~\eqref{eq:PF-Ritz-loss}, with a conventional data-driven supervised neural method.
The data-driven training is done with solutions given by the SAV method and the loss function is the difference between the neural-operator prediction and the SAV solution. The dataset is the same as for the unsupervised training. For both methods we use the same PR-RDNO architecture.

We further consider two different neural architectures for DeepRitzSplit.
\begin{itemize}
	\item 
    A PR-RDNO, as shown in Fig.~\ref{fig:RDNO_prescribed} and given by Eq.~\eqref{eq:aniso_scheme}. For the diffusion block a UNet with 4 downscaling and upscaling levels, 8 hidden channels, 3 width, periodic padding and hyperbolic tangent activation function is used. {The prescribed reaction term computes 
      \[ 
        \phi_n - \frac{\mathrm{dt}}{\tau} \nabla E_2(\phi_n, U_n).
        \]
    }
	\item
	A UNet with the same structure as in the diffusion block of the PR-RDNO described above.
\end{itemize}
The PR-RDNO architecture offers advantages in efficient training and improved interpretability of predictions, as it aligns with the scheme of Eq.~\eqref{eq:prescribed_react}, whereas the general UNet acts more as a black-box regression model.

All the training and evaluations were performed using a RTX6000 ADA GPU. The early stopping is applied once the loss fluctuates by less than 3\textperthousand.

\paragraph{Evaluation metrics}
To assess the accuracy of DeepRitzSplit approximations of the dendritic growth model, we evaluate the prediction of the neural operator by computing the evolution of a single grain. A solid particle is placed in the center of the unit square $[0, 1]^2$ at solidification temperature $U=0$ and with liquid at uniform temperature $\kappa<0$ (initial condition of Eq.~\eqref{eq:dendrite_initialcondition} for $n=1$ and $c_1 = (0.5, 0.5)$).

The reference to which we compare is the analytical solution for the growth speed and the radius of curvature of the dendrite tip, given by the theory of dendritic growth~\cite{Barbieri1989,dantzig2016solidification}. This theory shows that the tip of an isothermal dendrite branch growing into an infinite liquid at a lower temperature, settles at a constant growth speed, $V\tip$, and constant shape. The steady-state tip shape is a parabola with a radius of curvature $\rho\tip$ at the vertex. The solution of the Stefan problem for such a parabola~\cite{Ivantsov1947a} shows that the tip speed and radius are related to the dimensionless temperature difference that drives its growth, $\kappa$, by 
\begin{equation}
	\label{eq:Ivantsov}
	\sqrt{\pi Pe} \exp(Pe) \, \mathrm{erfc}(\sqrt{Pe}) + \kappa = 0 \,, 
\end{equation}
where $Pe = \rho\tip V\tip/(2 D)$ is the tip Péclet number, $D$ is the thermal diffusivity, and $\mathrm{erfc}$ is the complementary error function. A tip growing sufficiently far from the domain boundaries can be accurately approximated by this model and we use $Pe$ given by Eq.~\eqref{eq:Ivantsov} as a benchmark.

For the comparisons, we also track the approximated perimeters, given by Eq.~\eqref{eq:perimeter}, the free energy, $\Ed$, given by Eq.~\eqref{eq:free-energy}, and the evolution of the solid phase fraction, defined by 
\[
  \int_\Omega \left(\frac{1+\phi^\epsilon}{2}\right) d\vx \,.
\]

\begin{figure}[H]
	\centering
	\includegraphics[width=\textwidth]{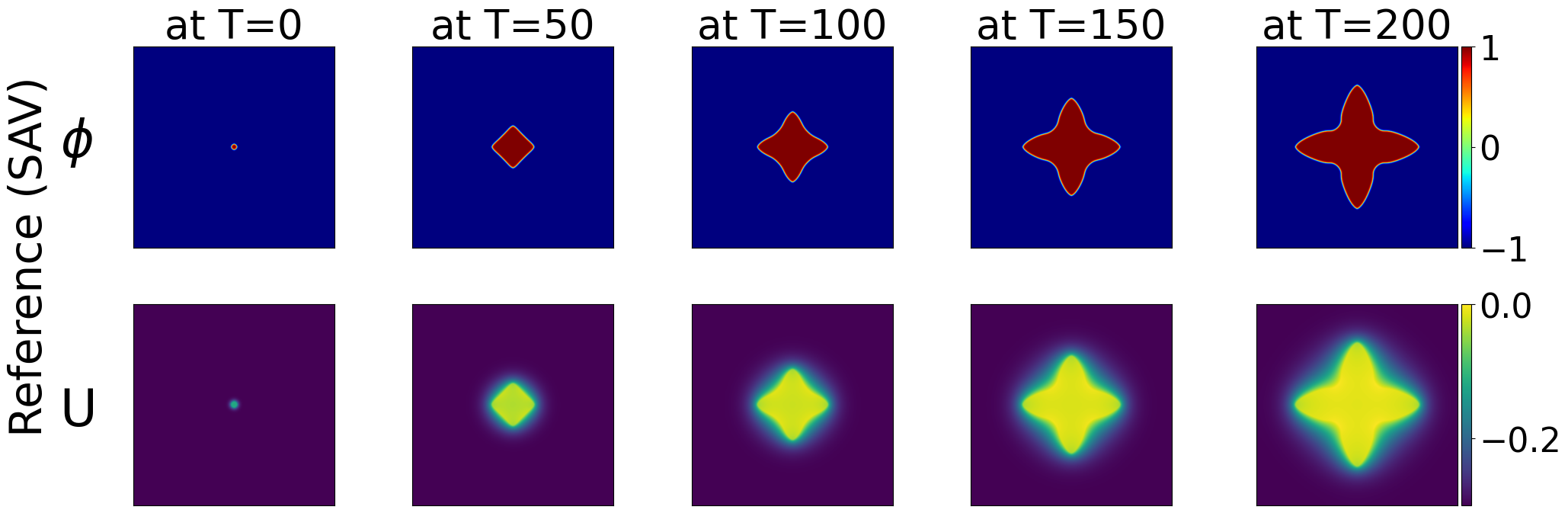}
	\includegraphics[width=\textwidth]{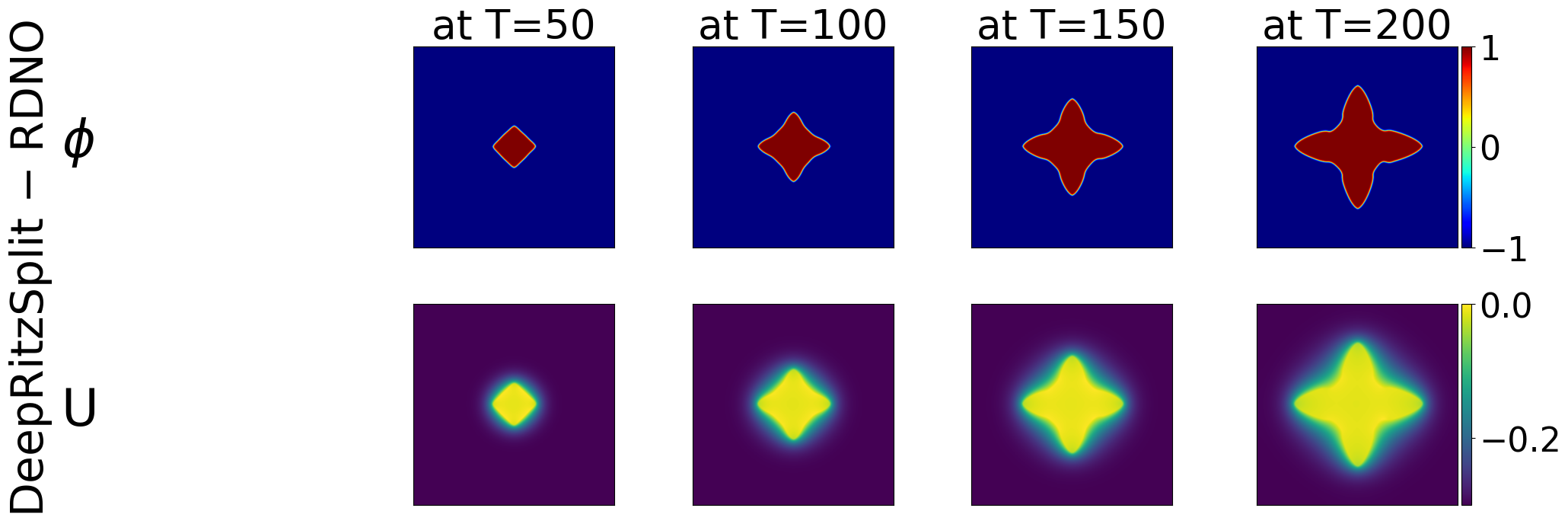}
	\includegraphics[width=\textwidth]{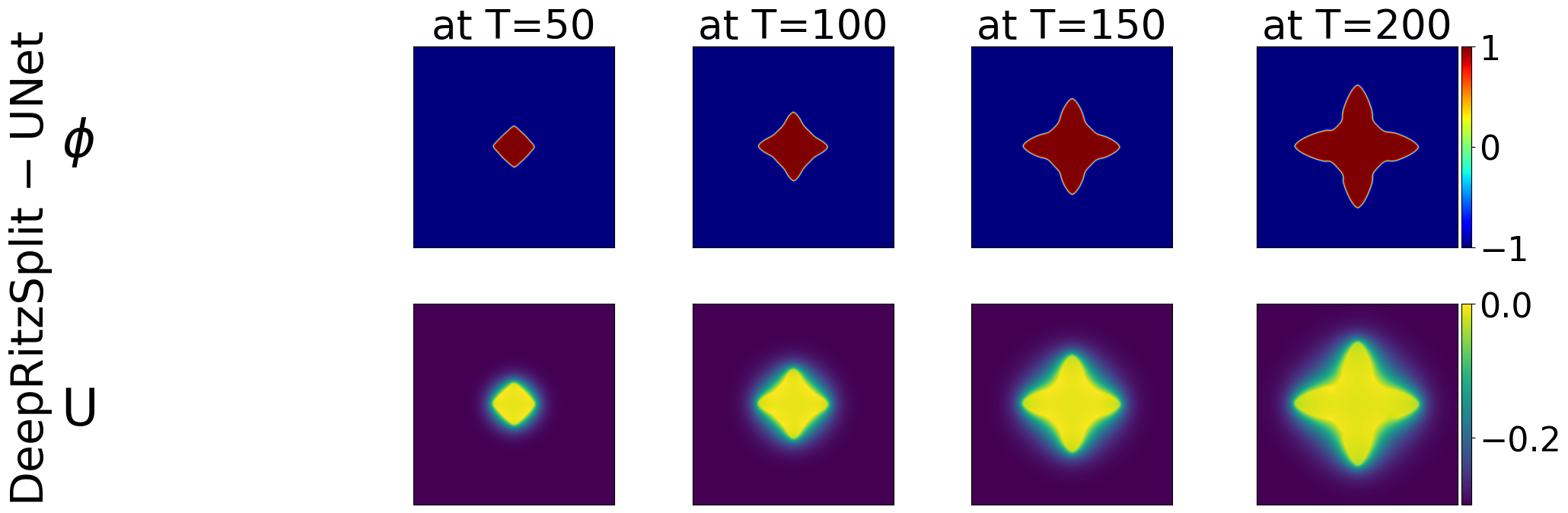}
	\includegraphics[width=\textwidth]{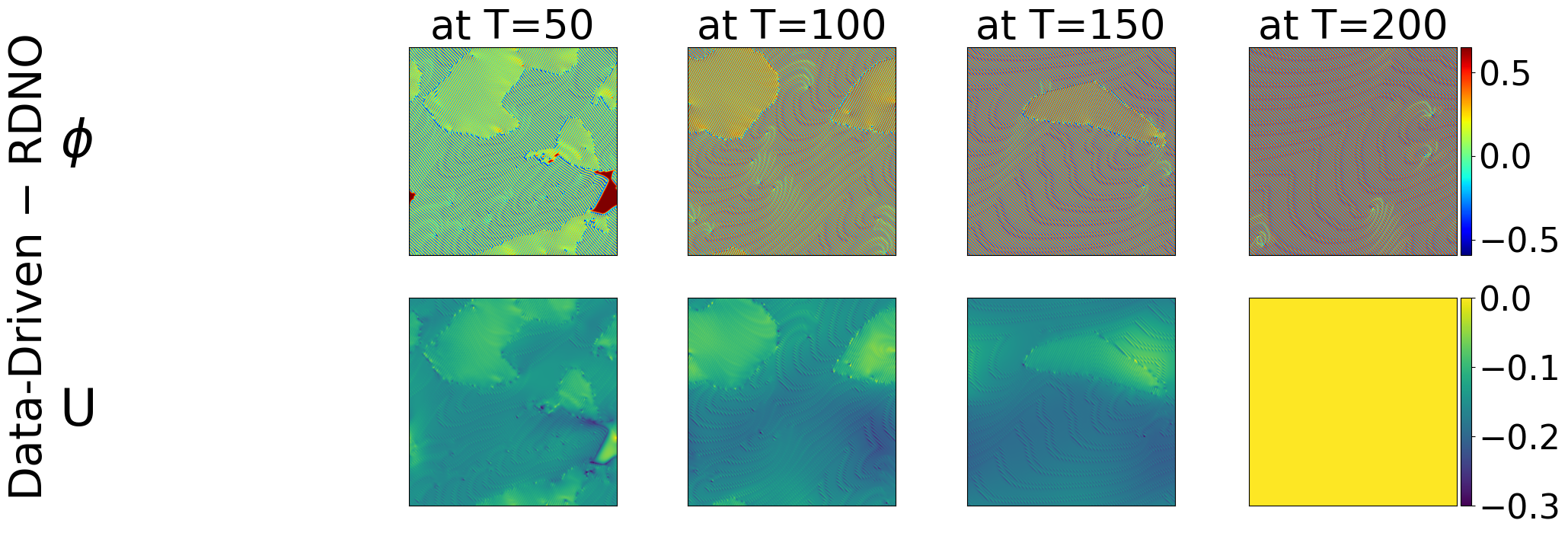}
  \caption{Single-grain dendritic growth at $\sigma=0.05$: SAV, DeepRitzSplit-based PR-RDNO (labelled RDNO in the figure), DeepRitzSplit-based UNet and supervised operator (top to bottom).}
	\label{fig:SAV_RDNO_UNet-one-gem.png}
\end{figure}

\begin{figure}[H]
  \centering
  \includegraphics[width=\textwidth]{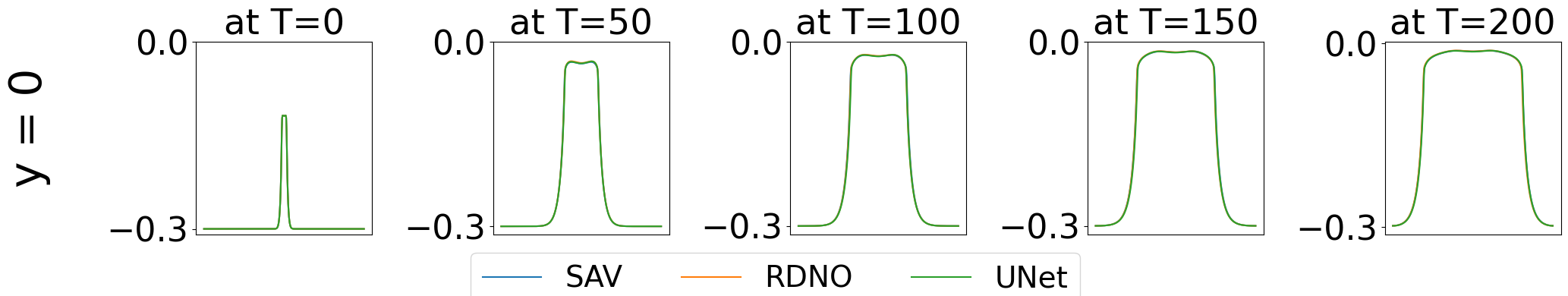}
  \includegraphics[width=\textwidth]{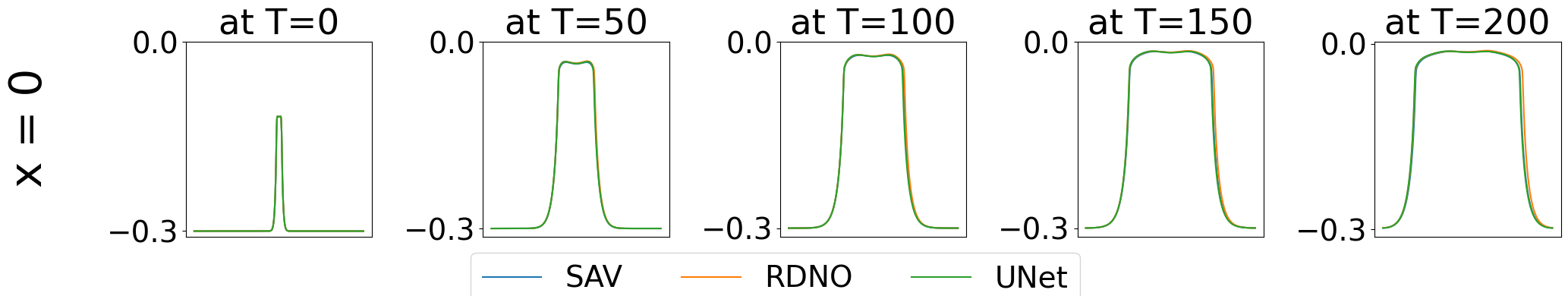}
  \includegraphics[width=\textwidth]{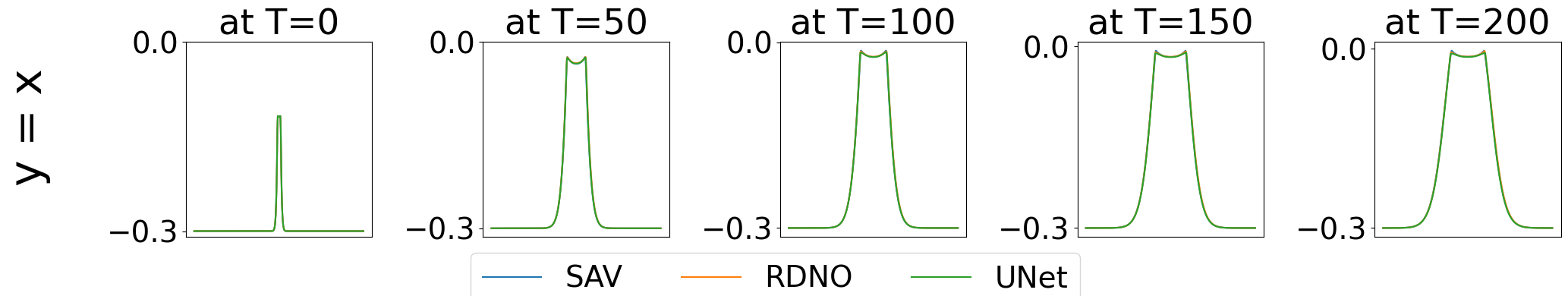}
  \caption{Comparison of dimensionless temperature distributions along the x-axis, y-axis, and the $y=x$ diagonal.}
  \label{fig:dimensionless_field_1d}
\end{figure}

\paragraph{Results} Comparisons of results for a single grain with $\sigma=0.05$ obtained with different neural operators are shown in Fig.~\ref{fig:SAV_RDNO_UNet-one-gem.png}. 
  Dimensionless temperature fields along three different axes ($y=0$, $x=0$ and $y=x$) are shown in Fig.~\ref{fig:dimensionless_field_1d}.
Note that the neural operator trained with a conventional data-driven approach completely fails to capture the solution. DeepRitzSplit with both 
PR-RDNO
and UNet performs relatively well in our test case, with 
PR-RDNO
having a somewhat better approximation of growth rate and better symmetry. Similar conclusions can be made for lower anisotropy, $\sigma=0.01$ (Fig.~\ref{fig:SAV_RDNO_UNet-one-gem-eps_m001.png} in Appendix~\ref{app:pf}).

The training times for the dendritic growth model, given in Table~\ref{tab:computation_times} are much longer than for the Allen-Cahn equation and require around 8 hours until the loss converges. The evaluation times are measured over 5000 iterations and with both neural architectures DeepRitzSplit is faster than the SAV numerical reference.

Relative errors of the evolution of solid phase fraction, perimeter and the free energy, at various time points, are given in Table~\ref{tab:pf_RDNO_UNET_rel_err}. All errors are within a few percent. The largest errors are in the perimeter and are due to somewhat more slender branches of the dendrite with the DeepRitzSplit - PR-RDNO at long times.

\begin{table}[H]
\caption{Training time for neural operators with DeepRitzSplit and comparison of evaluation times to the SAV numerical reference over 5000 iterations.}
\label{tab:computation_times}
\centering
\begin{tabular}{lccc}
	\hline
	        & {\bf Trainable } & \multicolumn{2}{c}{\bf Time} \\
	        & {\bf Parameters} & Train(hrs)  &    Eval.(s)    \\ \hline
	PR-RDNO &      589326      & $\approx$ 8 &      184       \\
	UNet    &      589326      & $\approx$ 8 &      183       \\
	SAV     &       n/a        &     n/a     &      269       \\ \hline
\end{tabular}
\end{table}

\begin{table}[H]
  \caption{Relative errors of DeepRitzSplit with the PR-RDNO and UNet architectures for 
  	the single grain configuration. Errors are given with respect to the SAV scheme reference solutions at different times and for both anisotropy strengths, $\sigma$.}
\label{tab:pf_RDNO_UNET_rel_err}
  \centering
  \textbf{DeepRitzSplit Relative Error, $\sigma = 0.01$}
\begin{tabular}{l|rr|rr|rr}
	\hline
	\multicolumn{1}{c|}{\bf } & \multicolumn{2}{c|}{\bf Solid fraction} & \multicolumn{2}{c|}{\bf Perimeter} & \multicolumn{2}{c}{\bf Energy} \\
	                          &     T=150 &                       T=450 &     T=150 &                  T=450 &     T=150 &              T=450 \\ \hline
	PR-RDNO                   &  $2.13\%$ &                    $1.80\%$ & $-0.77\%$ &              $-2.12\%$ & $-0.96\%$ &          $-2.37\%$ \\
	UNet                      & $-0.65\%$ &                   $-1.46\%$ & $-2.12\%$ &              $-4.59\%$ & $-0.48\%$ &           $0.71\%$ \\ \hline
\end{tabular}

\textbf{DeepRitzSplit Relative Error, $\sigma = 0.05$}
\begin{tabular}{l|rr|rr|rr}
	\hline
	\multicolumn{1}{c|}{\bf } & \multicolumn{2}{c|}{\bf Solid fraction} & \multicolumn{2}{c|}{\bf Perimeter} & \multicolumn{2}{c}{\bf Energy} \\
	                          &     T=100 &                       T=200 &     T=100 &                  T=200 &     T=100 &              T=200 \\ \hline
	PR-RDNO                   &  $2.03\%$ &                    $1.89\%$ &  $2.29\%$ &               $5.92\%$ & $-0.81\%$ &          $-2.17\%$ \\
	UNet                      & $-2.08\%$ &                   $-1.49\%$ & $-3.11\%$ &              $-1.15\%$ & $-0.37\%$ &          $-0.30\%$ \\ \hline
\end{tabular}
\end{table}

\begin{figure}[H]
	\centering
	\textbf{Anisotropic strength $\sigma=0.01$}
	\includegraphics[width=\textwidth]{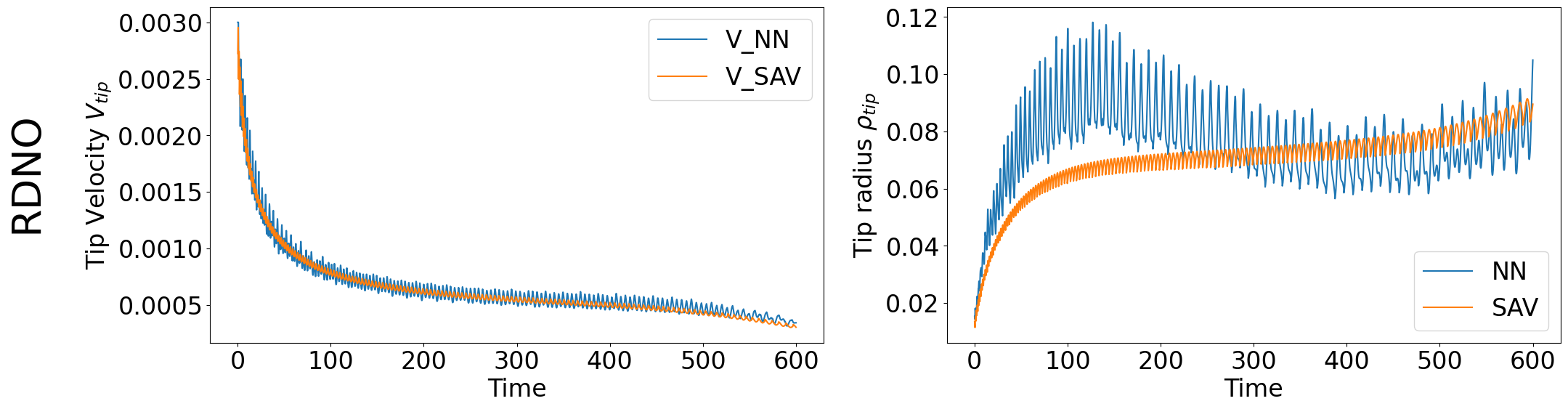}
	\includegraphics[width=\textwidth]{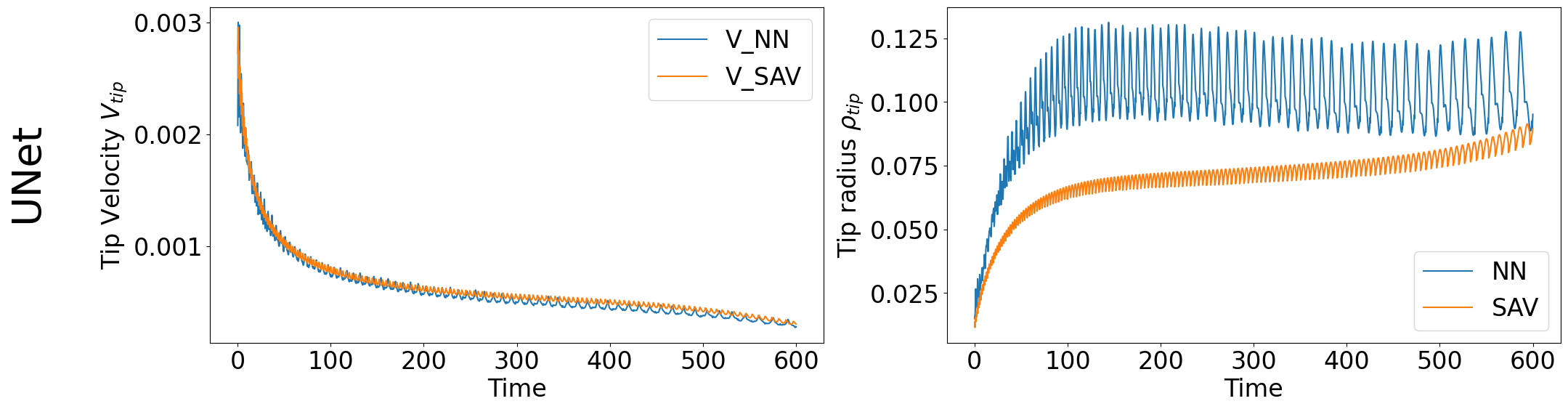}
	\\
	\textbf{Anisotropic strength $\sigma=0.05$}
	\includegraphics[width=\textwidth]{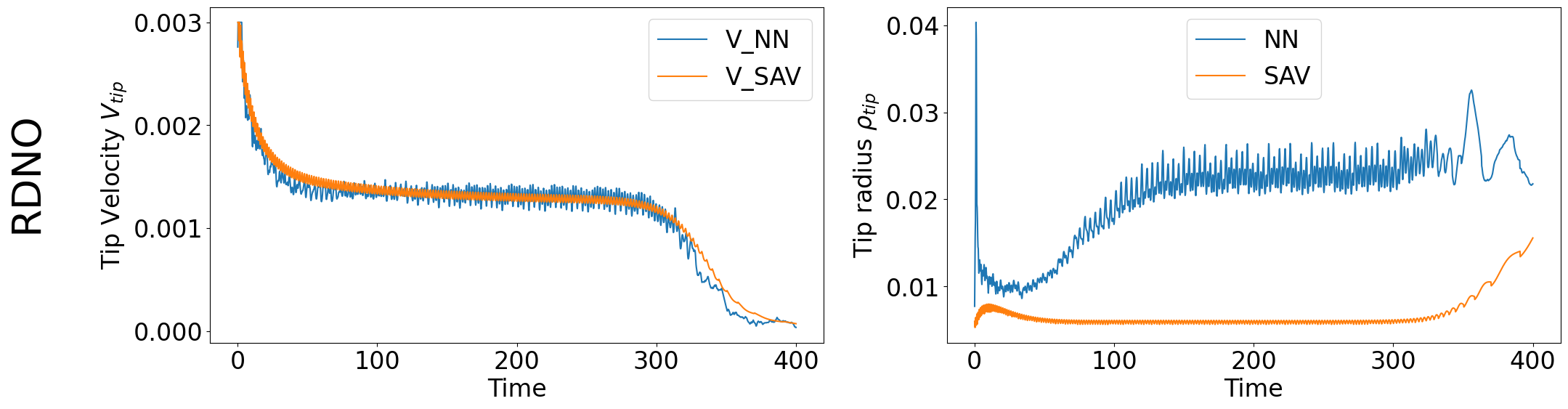}
	\\
	\includegraphics[width=\textwidth]{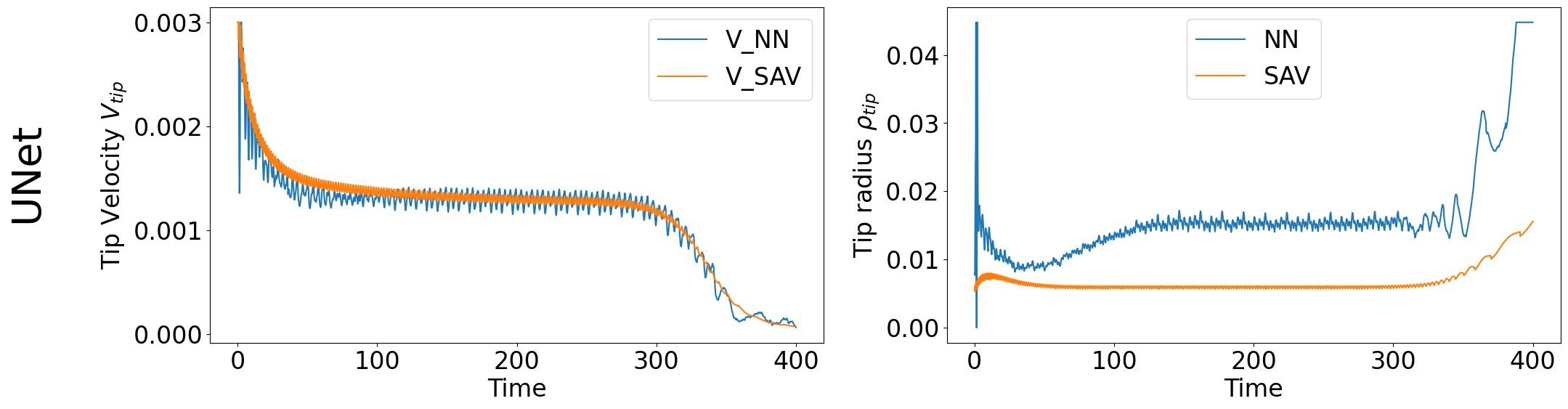}
	\caption{The evolution of tip velocity $V\tip$, and tip radius $\rho\tip$  with two different neural operator architectures (PR-RDNO and UNet) to the numerical reference (SAV) for simulations with a single grain. Note that the  growth slows down when the tip approaches the domain boundary.}
	\label{fig:RDNO_UNet_peclet}
\end{figure}

Next, we check whether these predictions are capable of capturing the underlying physics by evaluating the tip velocities and tip radius in the steady state. The evolutions of $V\tip$ and $\rho\tip$ are given in Fig.~\ref{fig:RDNO_UNet_peclet}. The DeepRitzSplit predictions show very good accuracy, closely following the velocity evolution predicted by the SAV numerical reference.
The differences of the methods are however very pronounced for the tip radius. The evolution strongly depends on the method (neural or 
numerical reference)
and on the neural architecture. 

The results on the steady-state dendrite tip operating state (velocity and radius), presented in Table~\ref{tab:UNet_RDNO_benchmark-peclet}, indicate that all approaches, the SAV scheme and both neural operators, give close tip velocity predictions but they also \textit{all} fail to properly characterize the dendrite tip radius.

\begin{table}[H]
  \caption{Velocity and radius of curvature of the dendrite tip during the steady\hyp state growth stage. Comparison of DeepRitzSplit with the PR-RDNO and UNet architectures to the numerical (SAV) and analytical references for simulations with a single grain. \added[id=CK]{$V_{tip}$ and $\rho_{tip}$ are computed as the average of 20 values around $T=400$ for $\sigma=0.01$ and around $T=200$ for $\sigma=0.05$, respectively.}}
\label{tab:UNet_RDNO_benchmark-peclet}
  \centering
  \textbf{Anisotropy strength $\sigma=0.01$}
\begin{tabular}{lccccrr}
	\hline
	                       &  \multicolumn{2}{c}{\textbf{Reference}}   & \multicolumn{2}{c}{\textbf{DeepRitzSplit}} & \multicolumn{2}{c}{\textbf{Relative error}} \\
	                       &     Analytical      &         SAV         &       PR-RDNO       &         UNet         &  PR-RDNO &                             UNet \\ \hline
	Tip speed, $V\tip$     &         n/a         & $5.45\cdot 10^{-4}$ & $5.57\cdot 10^{-4}$ & $4.94\cdot 10^{-4}$  &  $2.2\%$ &                          $9.3\%$ \\ 
	Tip radius, $\rho\tip$ &         n/a         & $7.20\cdot 10^{-2}$ & $7.93\cdot 10^{-2}$ & $6.47\cdot 10^{-3}$  & $10.2\%$ &                         $39.4\%$ \\ 
	Péclet number, $Pe$    & $4.48\cdot 10^{-2}$ & $3.14\cdot 10^{-1}$ & $3.53\cdot 10^{-1}$ & $2.56\cdot 10^{-1}$  & $12.6\%$ &                         $18.5\%$ \\ \hline
\end{tabular}

\textbf{Anisotropy strength $\sigma=0.05$}
\begin{tabular}{lccccrr}
	\hline
	                       &  \multicolumn{2}{c}{\textbf{Reference}}   & \multicolumn{2}{c}{\textbf{DeepRitzSplit}} & \multicolumn{2}{c}{\textbf{Relative error}} \\
	                       &     Analytical      &         SAV         &       PR-RDNO       &         UNet         &   PR-RDNO &                            UNet \\ \hline
	Tip speed, $V\tip$     &         n/a         & $1.29\cdot 10^{-3}$ & $1.32\cdot 10^{-3}$ & $1.30\cdot 10^{-3}$  & $-0.60\%$ &                       $-1.24\%$ \\
	Tip radius, $\rho\tip$ &         n/a         & $5.92\cdot 10^{-3}$ & $2.25\cdot 10^{-2}$ & $1.52\cdot 10^{-2}$  & $280.2\%$ &                       $157.0\%$ \\
	Péclet number, $Pe$    & $4.48\cdot 10^{-2}$ & $6.12\cdot 10^{-2}$ & $2.38\cdot 10^{-1}$ & $1.58\cdot 10^{-2}$  & $283.9\%$ &                       $158.7\%$ \\ \hline
\end{tabular}
\end{table}

\begin{figure}[H]
	\centering
	\includegraphics[width=\linewidth]{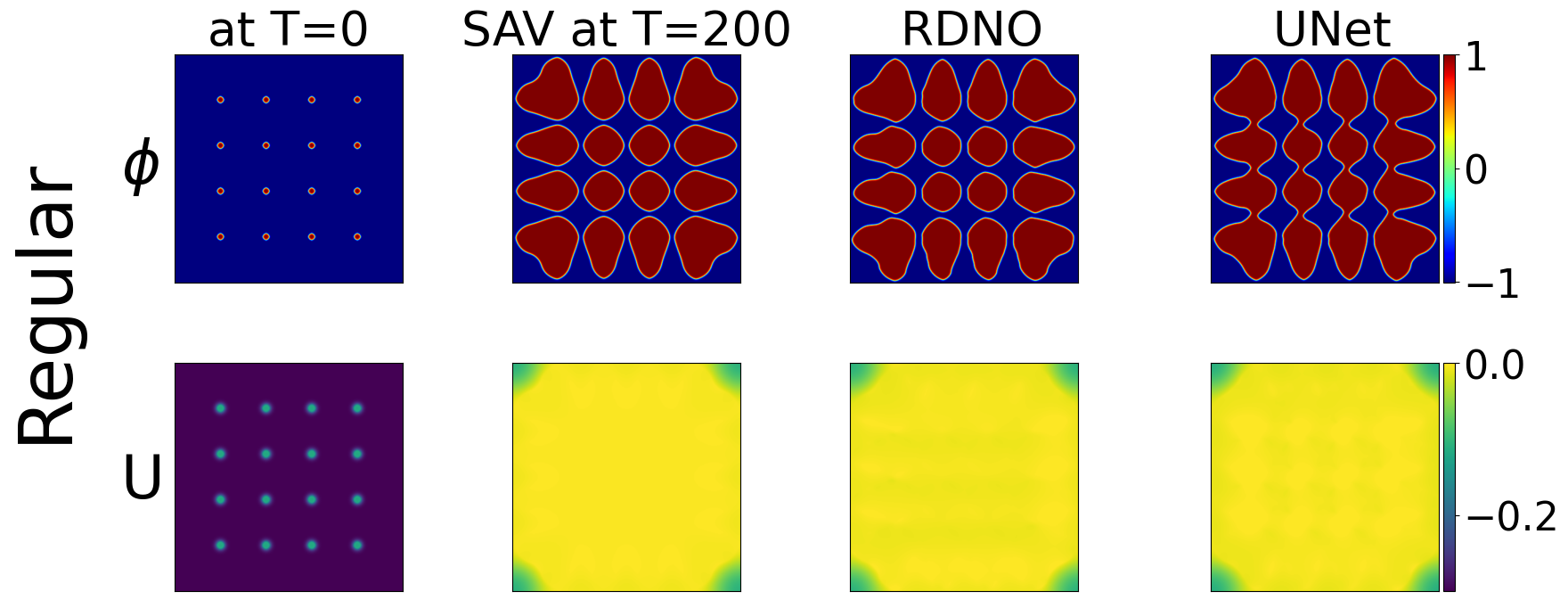}
	\includegraphics[width=\linewidth]{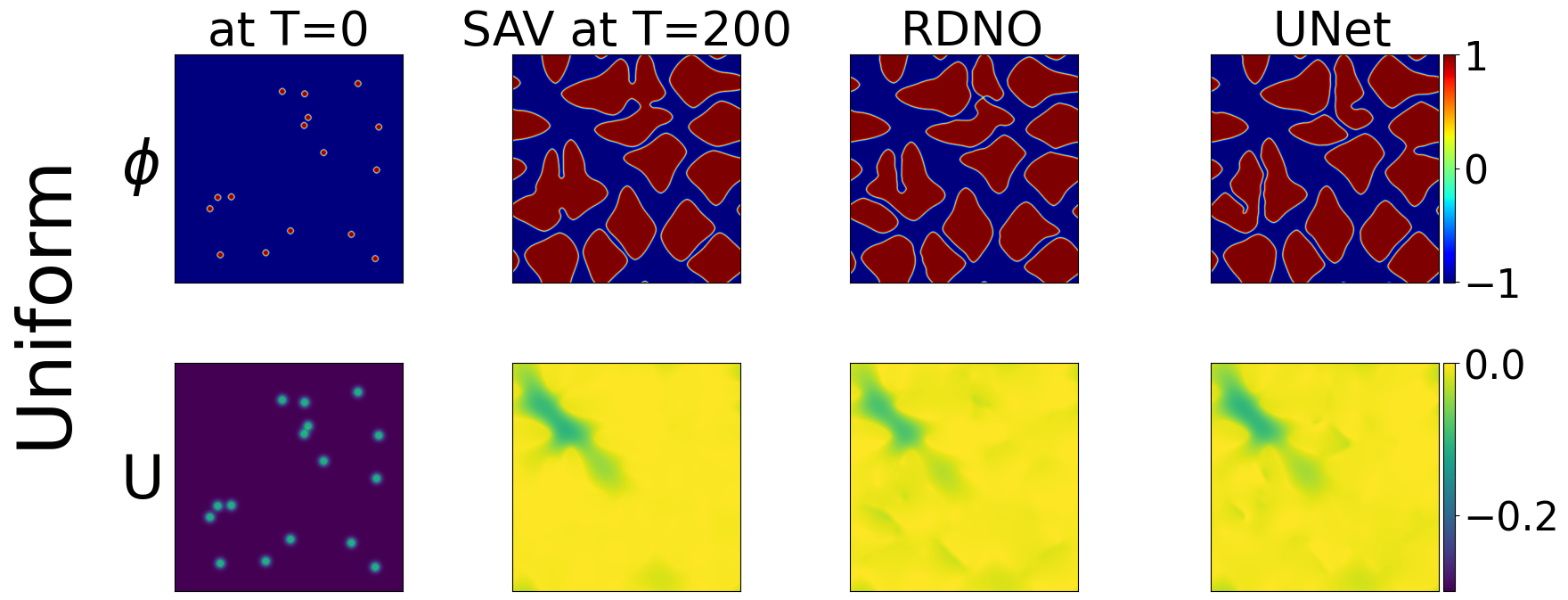}
	\includegraphics[width=\linewidth]{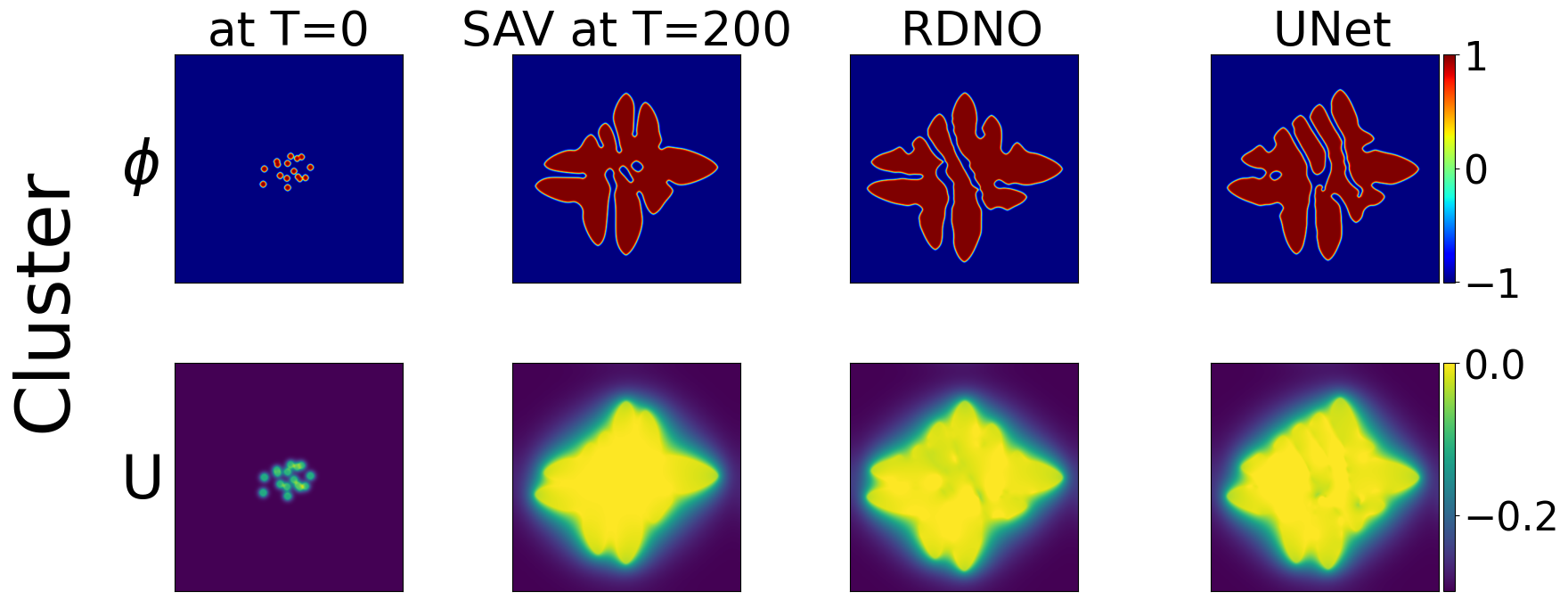}
  \caption{Sixteen-grain extrapolation at $\sigma=0.05$ for regular (top), uniform (mid) and clustered (bottom) arrangements. Each configuration compares SAV, PR-RDNO (labelled RDNO in the figure) and UNet.}
	\label{fig:multi-grains-strong}
\end{figure}

\begin{figure}[H]
	\centering
	\includegraphics[width=\linewidth]{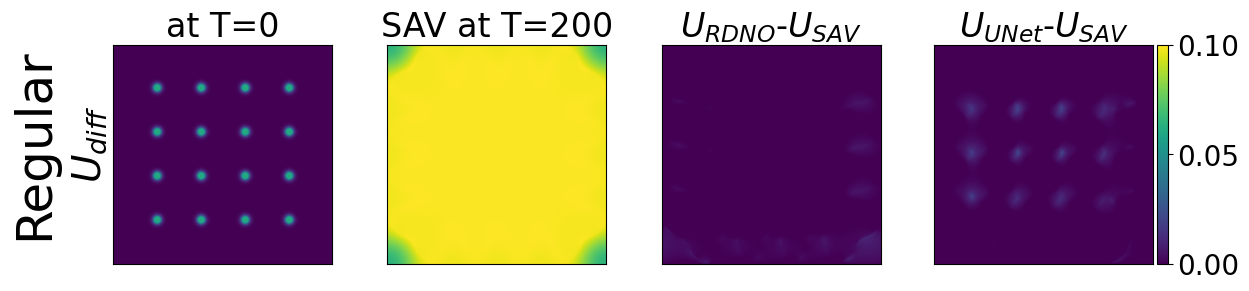}
	\includegraphics[width=\linewidth]{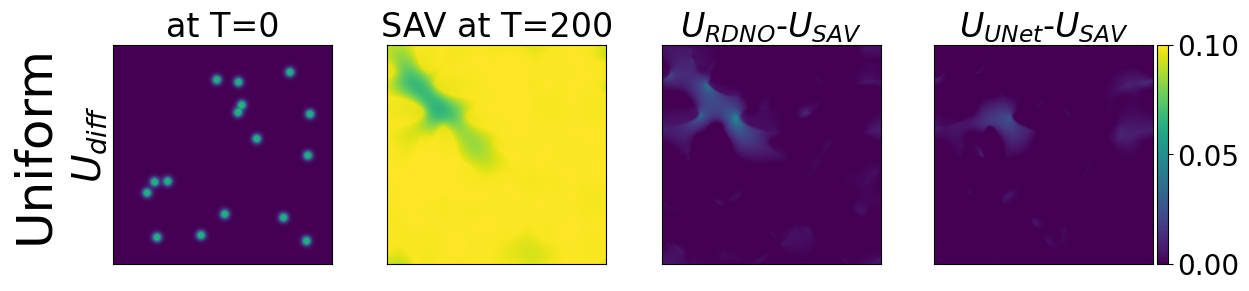}
	\includegraphics[width=\linewidth]{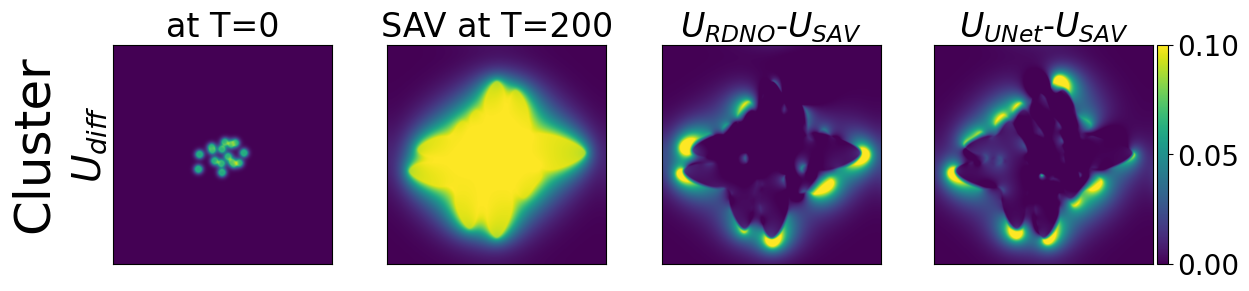}
  \caption{ 
Sixteen-grain extrapolation at $\sigma=0.05$ for regular (top), uniform (middle), and clustered (bottom) arrangements. Each configuration shows the dimensionless temperature field $U$ computed using SAV at T = 0 and T = 200, followed by the differences between the SAV results and those predicted by PR-RDNO (labelled RDNO in the figure) and UNet.
  }
  \label{fig:multi-Udiff}
\end{figure}

\paragraph{Application to multi-grain simulations}
For out-of-distribution evaluation, we test the trained neural operators on simulations with a larger number of grains than those used in the training set. 
This is particularly interesting for applications, since it shows if training can be done on a limited number of configurations with a small number of grains and the trained neural operator can then be used to simulate a wide variety of grain ensembles and spatial arrangements. Simulations with 16 grains across three different arrangements (regular, uniform, and cluster) are provided in Fig.~\ref{fig:multi-grains-strong}. 
\added[id=CK]{
The differences in dimensionless temperature for PR-RDNO and UNet relative to the SAV method are shown in Fig.~\ref{fig:multi-Udiff}, while the corresponding maximum temperatures for each configuration are summarized in Table~\ref{tab:U_fields_multi_max}.
}

\replaced[id=CK]{
One can observe that the PR-RDNO predictions align more closely with the SAV scheme in terms of dendrite morphology and the distribution of the dimensionless temperature field. While all methods yield results that generally remain consistent with the SAV reference, discrepancies are observed in local morphological features and temperature overshoot, most remarkably in the merging of grains. However, the underlying reason for this requires further investigation.
}
{
One can observe that the PR-RDNO predictions align more closely with the SAV scheme in terms of dendrite morphology. Overall, all methods give very similar results, close to the reference (SAV). The main differences are observed in local shape features {and temperature overshot}, most remarkably in the merging of grains. However, the underlying reason for this requires further investigation.
}

\begin{table}[H]
\caption{Maximum dimensionless temperatures for multi-grain simulations across all methods and spatial arrangements for the anisotropic strength $\sigma=0.05$.}
 \label{tab:U_fields_multi_max}
\centering
\begin{tabular}{lccc}
  \hline
  \textbf{Spatial}  & \multicolumn{3}{c}{\bf Evaluation methods}
  \\ 
  \textbf{Arrangements} & SAV & PR-RDNO & UNet 
  \\ \hline 
  Regular & $0$ & $7.37\cdot 10^{-3}$ & $3.16\cdot 10^{-1}$  
  \\ 
  Uniform & $3.06 \cdot 10^{-2}$ &  $4.47\cdot 10^{-1}$ & $5.57\cdot 10^{-1}$
  \\
  Cluster & $2.85\cdot 10^{-2}$ & $2.93\cdot 10^{-1}$ & $5.52\cdot 10^{-1}$
  \\ \hline
\end{tabular}
\end{table}

\paragraph{Refinement via Splitting scheme residuals: DeepRitzSplit as a Preconditioner}

Since the tip velocity is a function of the first-order derivative and the tip radius is a function of the second-order derivative of the phase solution $\phi$, 
the DeepRitzSplit method approximates well the tip velocity but not the tip radius. This discrepancy is likely due to the fact that the variational formulation used in the loss function only accounts for first-order terms.

Therefore, one could restart to train the neural operator with the splitting scheme residuals of the splitting scheme \eqref{eq:aniso_scheme} after the DeepRitzSplit training. Namely, we define the loss function by
\begin{equation}
  \mathcal{L}_{scheme}(\phi,\psi, U):=
\left \| \tau \frac{\phi- \psi}{dt}  +\nabla_\phi E_1(\phi, U) + \nabla_\phi E_2(\psi, U) \right \|_2^2
\label{eq:loss_scheme}
\end{equation}

\begin{figure}[H]
  \centering
  \textbf{Anisotropic strength $\sigma=0.01$}
\includegraphics[width=\textwidth]{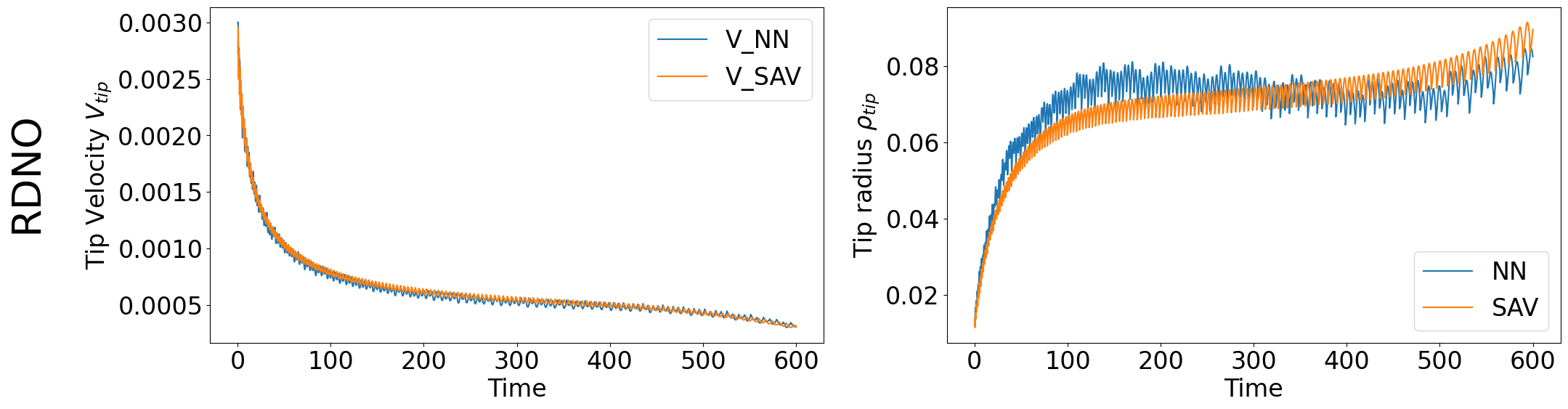}
\includegraphics[width=\textwidth]{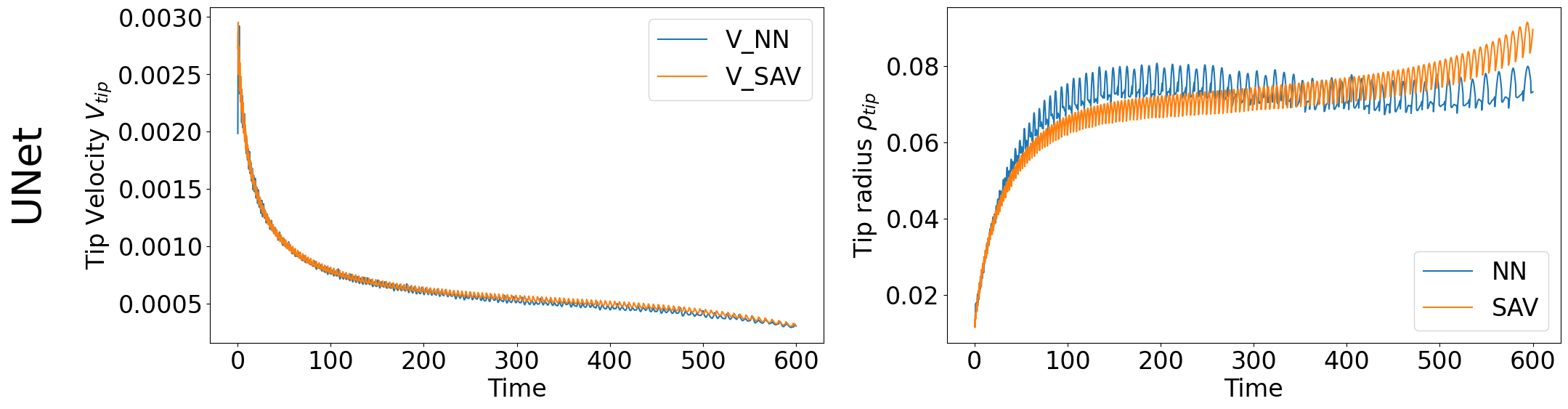}

  \textbf{Anisotropic strength $\sigma=0.05$}
\includegraphics[width=\textwidth]{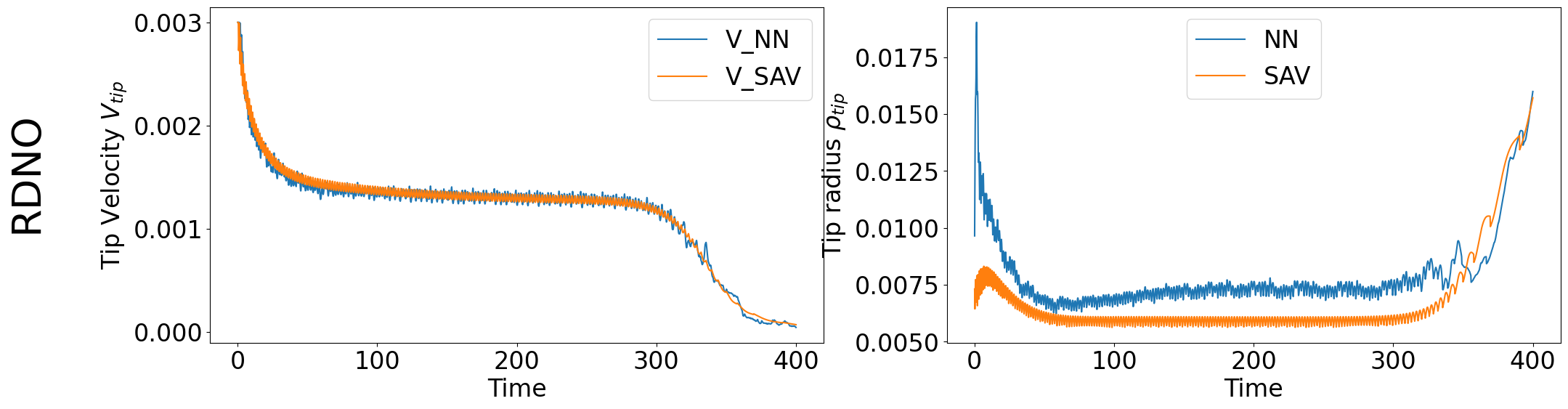}
\\
\includegraphics[width=\textwidth]{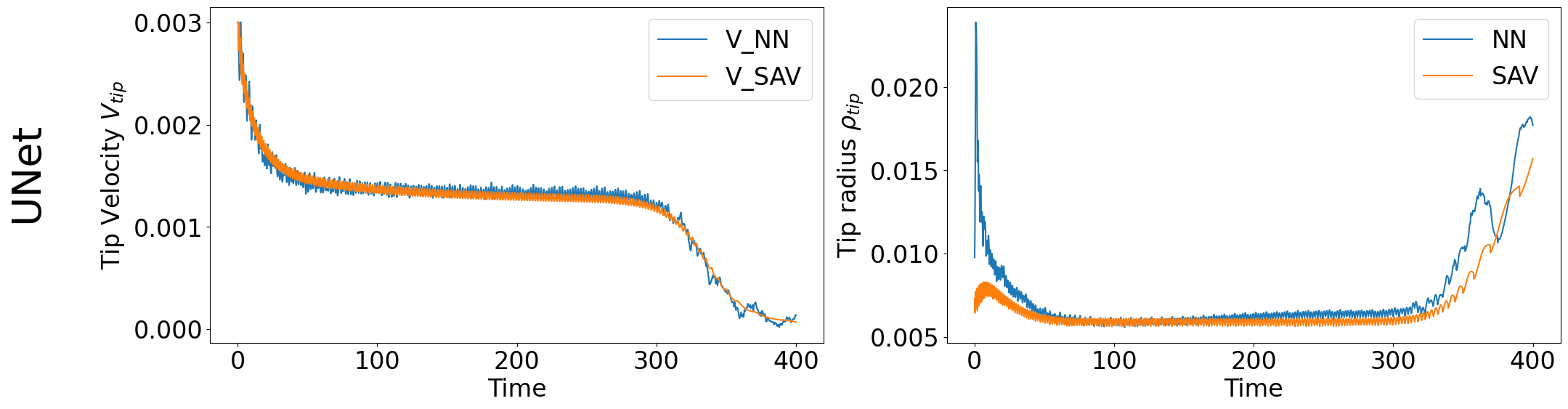}
\caption{Evolution of tip velocity $V\tip$ and tip radius $\rho\tip$ using PR-RDNO and UNet trained with DeepRitzSplit and scheme residuals~\eqref{eq:loss_scheme}, compared to the numerical reference (SAV) for simulations with a single grain. Note that the dendritic growth slow down when it starts to fill the whole domain.}
 \label{fig:RDNO_UNet_peclet_retrained}
\end{figure}

However, due to the high nonlinearity of the loss function~\eqref{eq:loss_scheme}, the training time required for convergence exceeds one week, compared to only a few hours for DeepRitzSplit.
Fortunately, as shown in Table~\ref{tab:UNet_RDNO_benchmark-peclet_retrained} and Fig.~\ref{fig:RDNO_UNet_peclet_retrained}, this refinement significantly improves the precision of the predicted tip radius.

\begin{table}[H]
  \caption{
Velocity and radius of curvature of the dendrite tip during the steady\hyp state growth stage.  
Here, neural operators are first pre-trained with DeepRitzSplit and subsequently refined using the scheme residuals~\eqref{eq:loss_scheme}.
{$V_{tip}$ and $\rho_{tip}$ are computed as the average of 20 values around $T=400$ for $\sigma=0.01$ and around $T=200$ for $\sigma=0.05$, respectively.}
  }
\label{tab:UNet_RDNO_benchmark-peclet_retrained}
  \centering
  \textbf{Anisotropic strength $\sigma=0.01$}

\begin{tabular}{lccccrr}
	\hline
	                       &       \multicolumn{2}{c}{Reference}        &     \multicolumn{2}{c}{DeepRitzSplit}     & \multicolumn{2}{c}{Relative error} \\
	                       &     Analytical      &         SAV          &       PR-RDNO       &        UNet         & PR-RDNO &                     UNet \\ \hline
	Tip speed, $V\tip$     &         n/a         & $5.45\times 10^{-4}$ & $5.28\cdot 10^{-4}$ & $5.28\cdot 10^{-4}$ & $3.4\%$ &                  $3.2\%$ \\ 
	Tip radius, $\rho\tip$ &         n/a         & $7.20\cdot 10^{-2}$  & $7.49\cdot 10^{-2}$ & $7.11\cdot 10^{-3}$ & $4.1\%$ &                  $1.2\%$ \\ 
	Péclet number, $Pe$    & $4.48\cdot 10^{-2}$ & $3.14\cdot 10^{-1}$  & $3.16\cdot 10^{-1}$ & $3.00\cdot 10^{-1}$ & $0.8\%$ &                  $4.3\%$ \\ \hline
\end{tabular}
\textbf{Anisotropic strength $\sigma=0.05$}
\begin{tabular}{lccccrr}
	\hline
	                       &       \multicolumn{2}{c}{Reference}       &     \multicolumn{2}{c}{DeepRitzSplit}     & \multicolumn{2}{c}{Relative error} \\
	                       &     Analytical      &         SAV         &       PR-RDNO       &        UNet         &  PR-RDNO &                    UNet \\ \hline
	Tip speed, $V\tip$     &         n/a         & $1.29\cdot 10^{-3}$ & $1.29\cdot 10^{-3}$ & $1.34\cdot 10^{-3}$ &  $0.3\%$ &                 $3.4\%$ \\
	Tip radius, $\rho\tip$ &         n/a         & $5.92\cdot 10^{-3}$ & $7.07\cdot 10^{-3}$ & $6.20\cdot 10^{-3}$ & $19.3\%$ &                 $4.7\%$ \\
	Péclet number, $Pe$    & $4.48\cdot 10^{-2}$ & $6.11\cdot 10^{-2}$ & $7.30\cdot 10^{-2}$ & $6.65\cdot 10^{-2}$ & $19.4\%$ &                 $8.8\%$ \\ \hline
\end{tabular}
\end{table}

We observe that if the neural operator is trained directly using the scheme residuals, the training hardly converges and the loss remains unstable. This observation motivates thus the use of DeepRitzSplit as a preconditioner, providing a ``warm-start'' that enables convergence within a reasonable timeframe.

\section{Conclusion and perspectives}
\label{sec:conclusion}

In this paper, we introduced DeepRitzSplit, a variational training framework for neural operators based on the minimization problem associated with an energy-stable convex--concave splitting scheme. The Allen--Cahn benchmark shows that, for identical neural-operator architectures, this unsupervised training strategy generalizes substantially better than conventional supervised data-driven training in the out-of-distribution configurations considered here. The trained operator approximates the time-advancement map of the splitting scheme while the loss retains its underlying dissipative variational structure.

We then applied the same framework to the substantially more complex phase-field model of anisotropic dendritic growth. The neural operators reproduce global quantities such as solid fraction and dendrite tip velocity with good accuracy and generalize from training configurations containing one to three grains to more complex multi-grain configurations. The results also expose an important limitation: quantities involving higher spatial derivatives, in particular the dendrite tip radius, are much less accurately captured by the initial DeepRitzSplit training. Refining a DeepRitzSplit-trained operator with the residual of the splitting scheme markedly improves this quantity, whereas training directly from the scheme residual is difficult to converge. DeepRitzSplit can therefore also be viewed as an effective physics-based initialization for a more accurate residual-based training stage.

The comparison of architectures further shows that the structure imposed on the neural operator matters. The RDNO used for Allen--Cahn mirrors the learned reaction and diffusion contributions of the splitting scheme, whereas the Prescribed-Reaction RDNO (PR-RDNO) used for dendritic growth evaluates the reaction contribution explicitly and learns the remaining implicit operator. This separation improves interpretability and provides a useful route for designing neural surrogates adapted to more complex gradient-flow models.

Further work should quantify more systematically the robustness of these operators with respect to anisotropy strength, undercooling, grain number and grain arrangement, and determine when the gain in evaluation time compensates for the training cost and the remaining approximation error.

\subsection*{Acknowledgements}
This work was supported by the Lorraine University of Excellence through the interdisciplinary project TransPINNO.)
Experiments presented in this paper were carried out using the Grid'5000 testbed, supported by a scientific interest group hosted by Inria and including CNRS, RENATER and several Universities as well as other organizations (see \url{https://www.grid5000.fr})

\printbibliography

\appendix
\section{Appendix}\label{sec:app}

\subsection{Neural Operator for the Allen-Cahn equation}\label{app:NO_AC}

%

\begin{figure}[H]
  \centering
    \includegraphics[width=.8\textwidth]{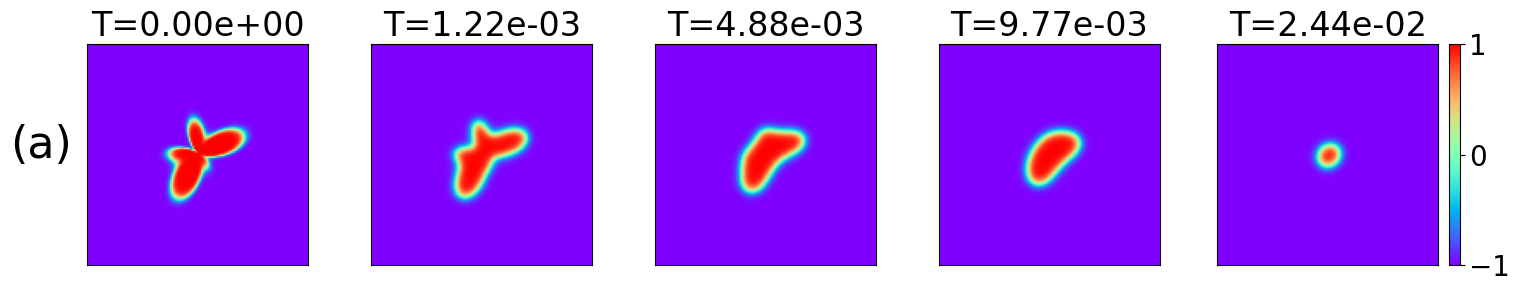}
    \includegraphics[width=.8\textwidth]{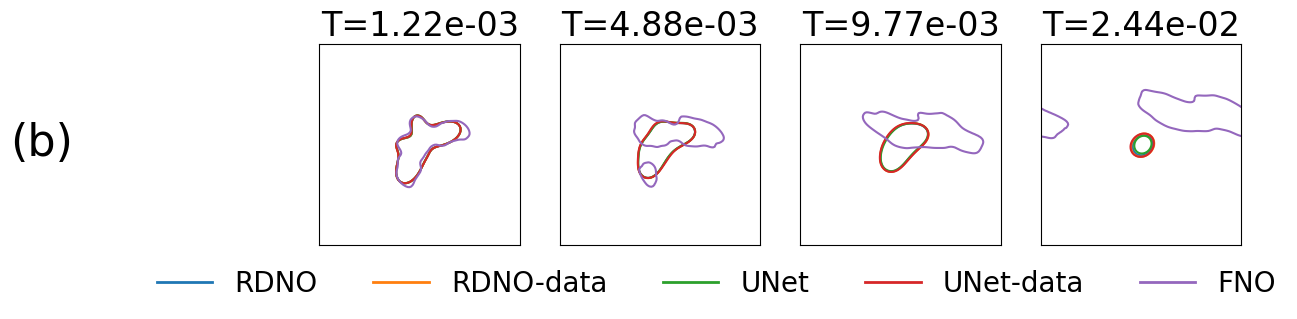}
    \includegraphics[width=.8\textwidth]{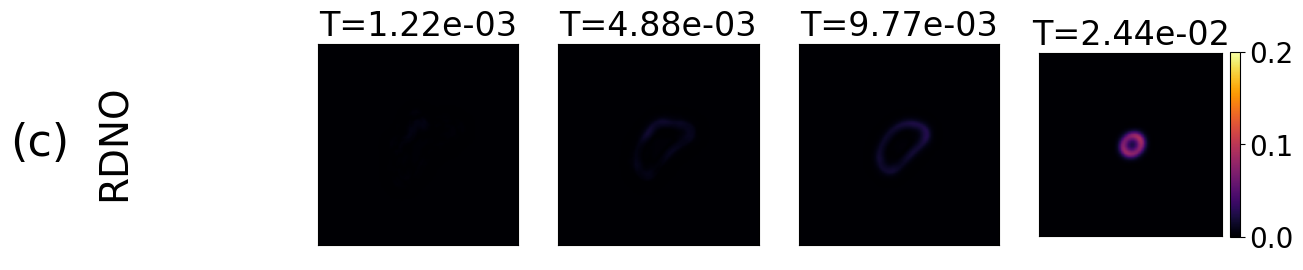}
    \includegraphics[width=.8\textwidth]{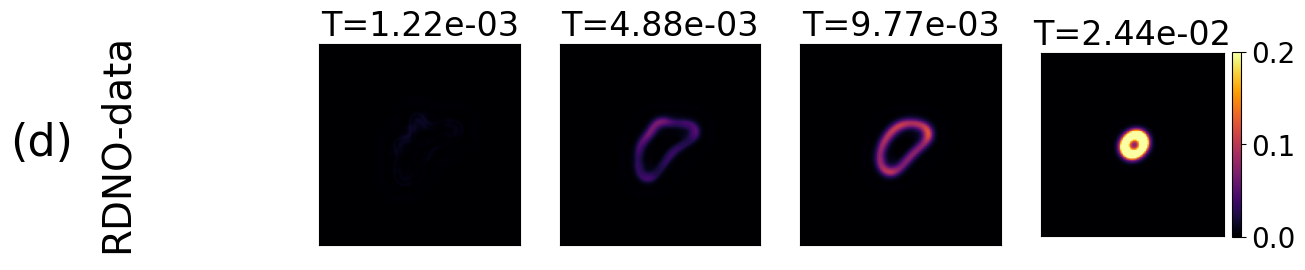}
    \includegraphics[width=.8\textwidth]{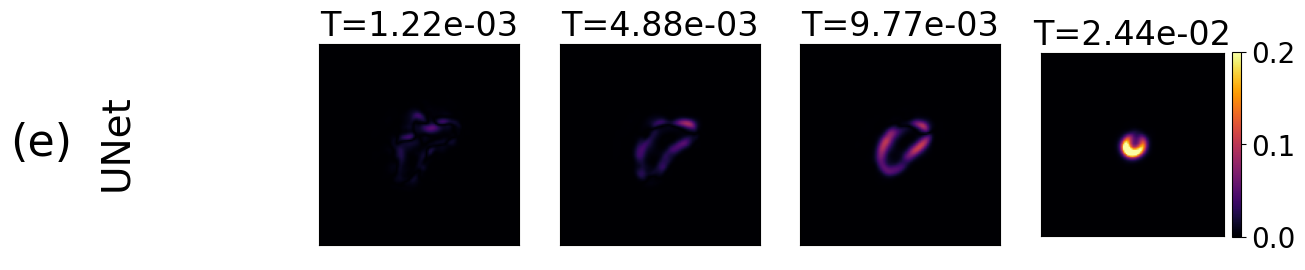}
    \includegraphics[width=.8\textwidth]{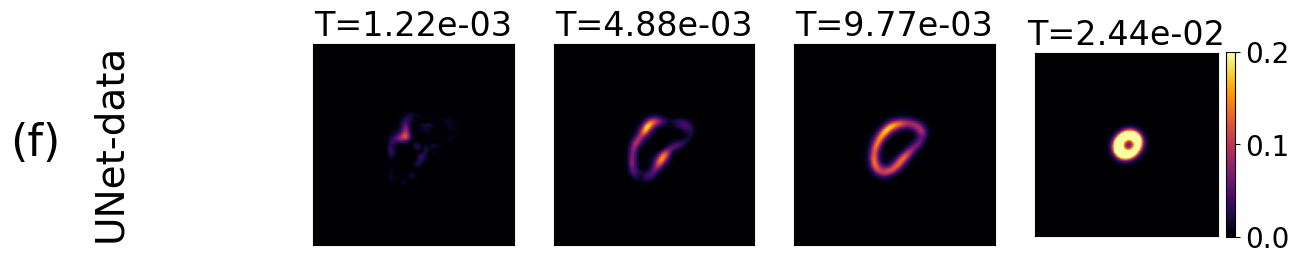}
    \includegraphics[width=.8\textwidth]{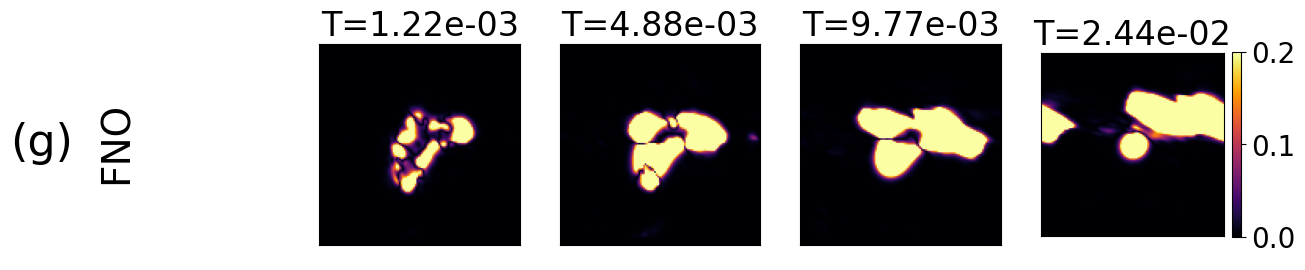}
    \includegraphics[width=.8\textwidth]{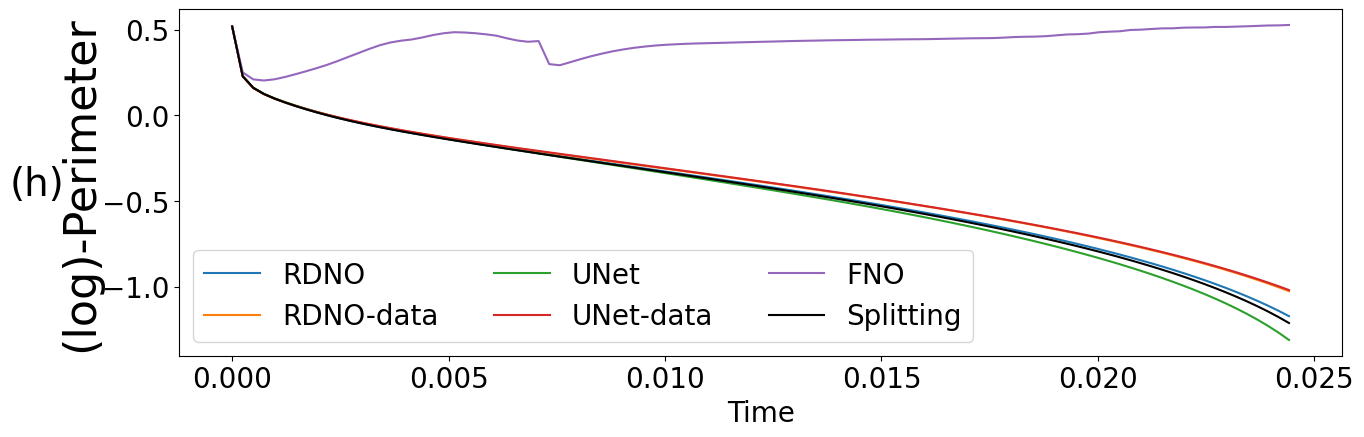}
  \caption{Predictions of the neural operator methods and error distribution on an in-distribution configuration}
  \label{fig:AC_IOD}
\end{figure}

\begin{figure}[H]
  \centering
    \includegraphics[width=.8\textwidth]{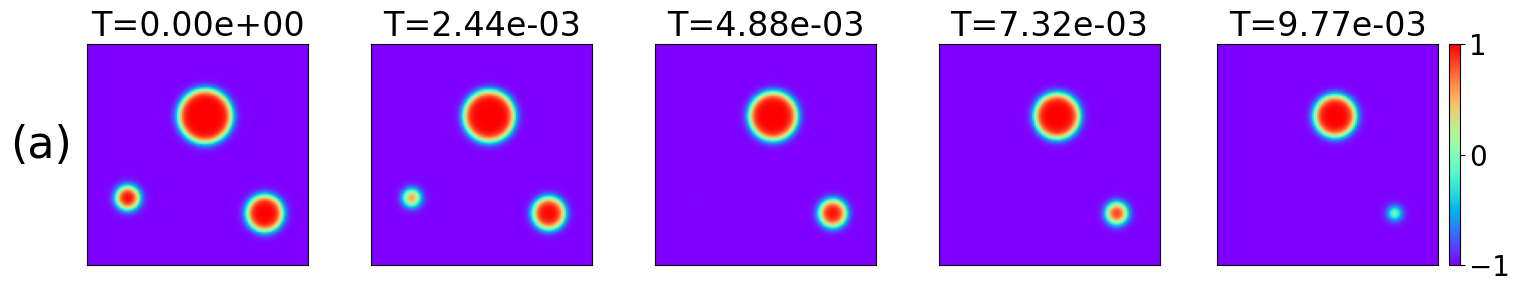}
    \includegraphics[width=.8\textwidth]{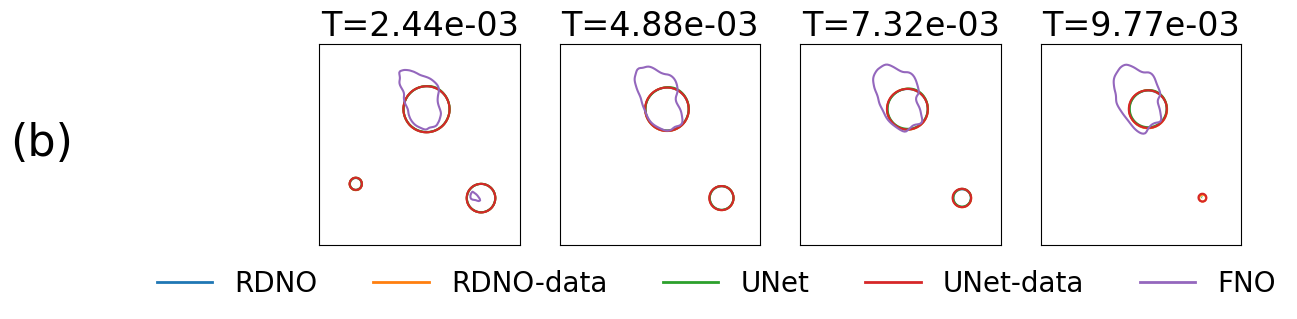}
    \includegraphics[width=.8\textwidth]{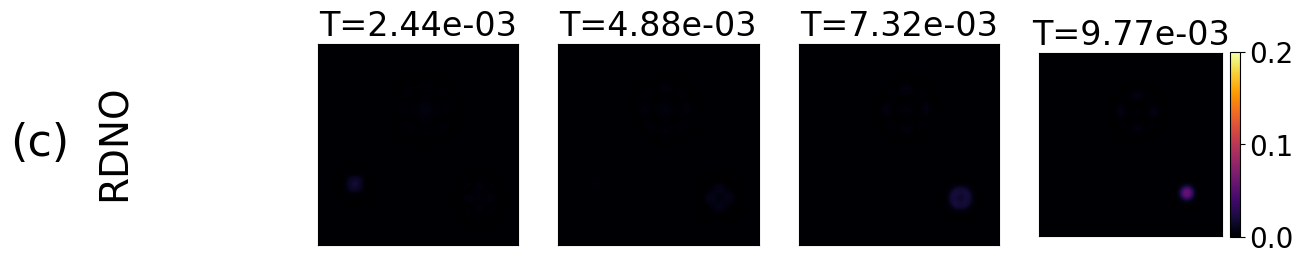}
    \includegraphics[width=.8\textwidth]{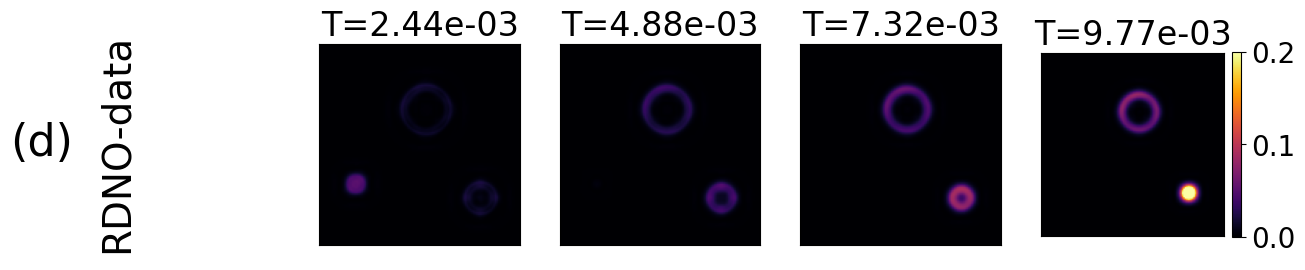}
    \includegraphics[width=.8\textwidth]{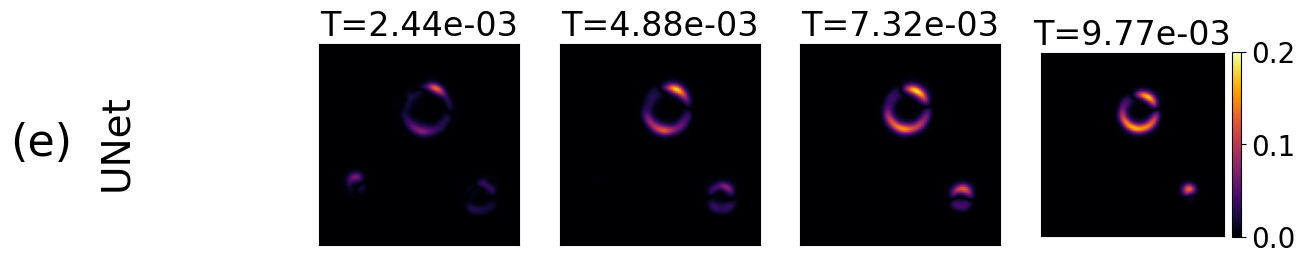}
    \includegraphics[width=.8\textwidth]{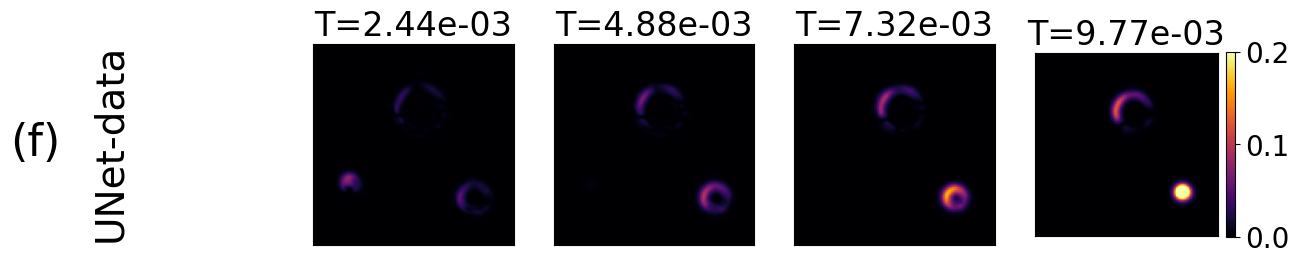}
    \includegraphics[width=.8\textwidth]{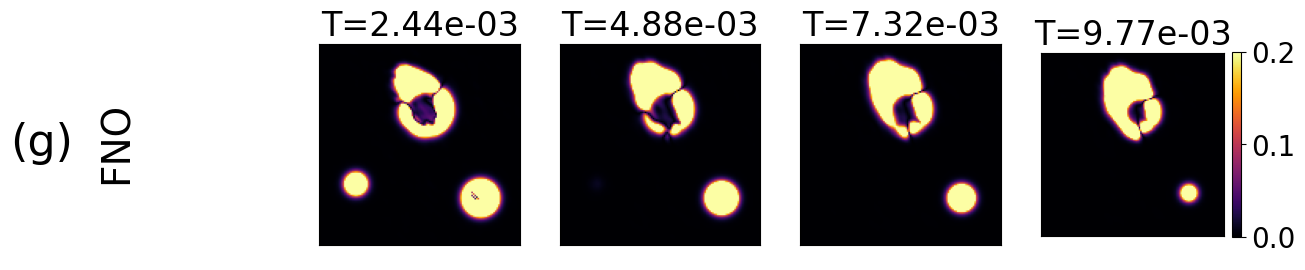}
    \includegraphics[width=.8\textwidth]{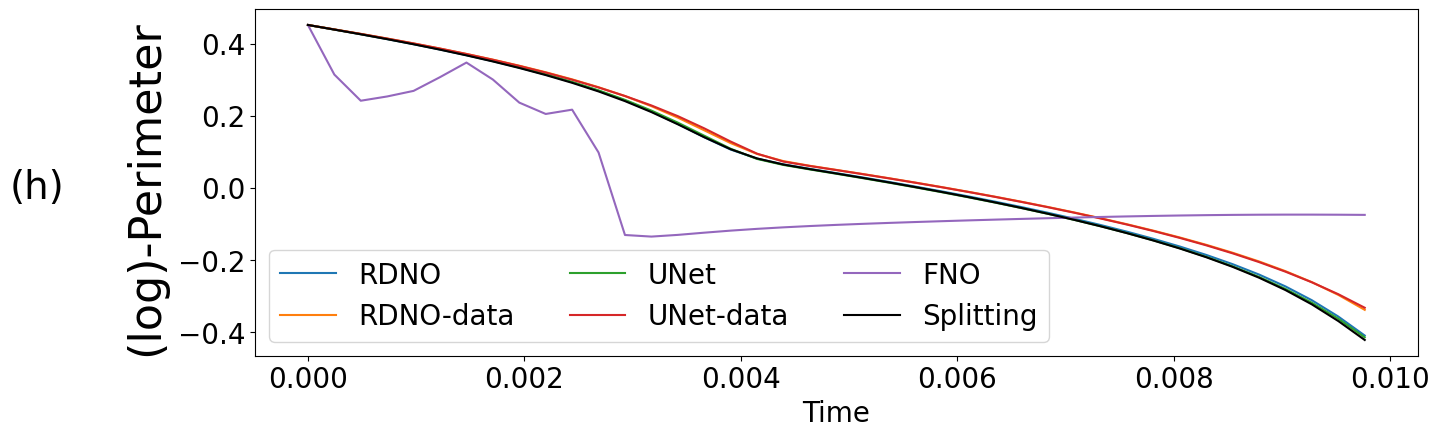}
  \caption{Predictions of the neural operator methods and error distribution on a multi-disk configuration}
  \label{fig:AC_multi}
\end{figure}

\subsection{SAV based method for anisotropic dendritic growth}
\label{app:sav}

The SAV method is a numerical scheme designed to handle gradient flows and energy-based PDEs while preserving their inherent energy dissipation properties. One of its main advantages is its \emph{unconditional energy stability}, which allows for larger time steps without sacrificing the physical accuracy of the simulation. Additionally, by reformulating complex nonlinear terms using an auxiliary scalar variable, the SAV method simplifies the numerical treatment and enables the use of efficient implicit time-stepping schemes. However, 
the introduction of auxiliary variables may increase the computational time, and the method's generalization to highly nonlinear systems can be nontrivial.

Initially proposed in \cite{li2022new} and later improved in \cite{guo2024efficient}, a SAV based splitting method was introduced for the anisotropic phase-field dendritic crystal growth model.
The core idea is to first consider the following decoupling of the free energy $E$:
$$ 
E(\phi , U) = \tilde E_1(\phi, U ) + \tilde E_2(\phi, U)
$$ 
where
$$ 
\tilde E_1(\phi, U) = \int_\Omega \left ( 
    \frac{1}{2} \alpha |\nabla \phi|^2+ \beta \frac{\phi^2}{2 \epsilon^2} 
\right )\, dx
$$
with $\alpha, \beta >0$
and
$$ 
\tilde E_2(\phi, U) = E(\phi, U) - E_1(\phi, U). 
$$
The corresponding splitting numerical scheme is therefore given by
\begin{equation}
  \left \{
\begin{aligned}
 \tau \frac{\phi_{n+1}- \phi_n}{dt} &=  -\nabla_\phi \tilde E_1(\phi_{n+1}, U_n) - \nabla_\phi \tilde E_2(\phi_n, U_n)
 \\ 
                                    &= \alpha \Delta \phi_{n+1} - \beta \epsilon^{-2} \phi_{n+1} - 
                                       \nabla_\phi \tilde E_2(\phi_n, U_n)
\\ 
 \frac{U_{n+1} - U_n}{dt} &= D \Delta U_{n+1} + K h'(\phi_{n+1}) \frac{\phi_{n+1} - \phi_n}{dt}
\end{aligned}
\right .
  \label{eq:simple_scheme}
\end{equation}

However, regardless of the numerical efficiency of solving Eq.~\eqref{eq:simple_scheme}, such a numerical scheme~\eqref{eq:simple_scheme} is \emph{not} unconditionally stable since the concavity of $\tilde E_2$ is not ensured in general. 
In order to guarantee energy stability, Guo et al. reformulate the scheme~\eqref{eq:simple_scheme} by adding a scalar auxiliary variable $q_n$ as an energy relaxation term, see Eqs.~\eqref {eq:SAV_scheme-1} and~\eqref{eq:SAV_scheme} below.
Despite requiring more computational time compared to Eq.~\eqref{eq:simple_scheme}, this approach enables the use of efficient implicit time-stepping schemes while ensuring energy dissipation as described in Eq.~\eqref{eq:energy_dissipation}.
Unlike traditional numerical methods that may struggle with stability or require small time steps for accuracy, the SAV method allows for larger time steps without compromising stability.

The SAV method introduced in \cite{guo2024efficient} consists of two steps. First, for $(\phi_n, U_n, q_n)$, we seek for a solution $(\bar{\phi}_{n+1}, \bar{U}_{n+1}, \bar{q}_{+1})$ of,
\begin{equation}
  \left \{
\begin{aligned}
 \tau \frac{\bar{\phi}_{n+1}- \phi_n}{\mathrm{dt}} &=  -\nabla_\phi \tilde E_1(\bar{\phi}_{n+1}, U_n) - \nabla_\phi \tilde E_2(\phi_n, U_n)
 \\ 
                                    &= \alpha \Delta \bar{\phi}_{n+1} - \beta \epsilon^{-2} \bar{\phi}_{n+1} - 
                                       \nabla_\phi \tilde E_2(\phi_n, U_n)
\\ 
 \frac{\bar{U}_{n+1} - U_n}{dt} &= D \Delta \bar{U}_{n+1} + K h'(\bar{\phi}_{n+1}) \frac{\bar{\phi}_{n+1} - \phi_n}{\mathrm{dt}}
 \\ 
 \frac{( \bar{q}_{n+1} - q_n)}{\mathrm{dt}}  &=  \xi_{n+1} dE_t(\bar{\phi}_{n+1}, \bar{U}_{n+1})
\end{aligned}
\right .
  \label{eq:SAV_scheme-1}
\end{equation}
with $q_0 = E(\phi_0, U_0)$ and $\xi_{n+1} = \frac{\bar{q}_{n+1}}{E(\bar{\phi}_{n+1}, \bar{U}_{n+1})}$.
 And then we set $\eta_{n+1} = 1- (1- \xi_{n+1})^2$ and 
\begin{equation}
  \left \{
\begin{aligned}
 \phi_{n+1} &= \eta_{n+1} \bar{\phi}_{n+1}
 \\
 U_{n+1} &= \eta_{n+1} \bar{U}_{n+1}
 \\ 
 q_{n+1} &= \zeta_{n+1} \bar{q}_{n+1} + (1 - \zeta_{n+1}) E(\phi_{n+1}, U_{n+1})
\end{aligned}
\right .
  \label{eq:SAV_scheme}
\end{equation}
where  $\zeta_{n+1} \in [0, 1] $ is chosen such that 
\begin{equation}
 \frac{( q_{n+1} - \bar{q}_{n+1})}{\mathrm{dt}} \leq  - \xi_{n+1} dE_t(\bar{\phi}_{n+1}, \bar{U}_{n+1}) \, .
  \label{eq:SAV-energy-decreasing}
\end{equation}
With the help of the auxiliary variable $q_n$, it is proven in \cite{guo2024efficient} that the numerical scheme \eqref{eq:SAV_scheme-1}-\eqref{eq:SAV_scheme} is energetically stable.

\subsection{Phase-Field model of dendritic growth }
\label{app:pf}

\textbf{Weak anisotropic strengh of surface energy $\sigma = 0.01$:}

\begin{figure}[H]
  \centering
  \includegraphics[width=.9\textwidth]{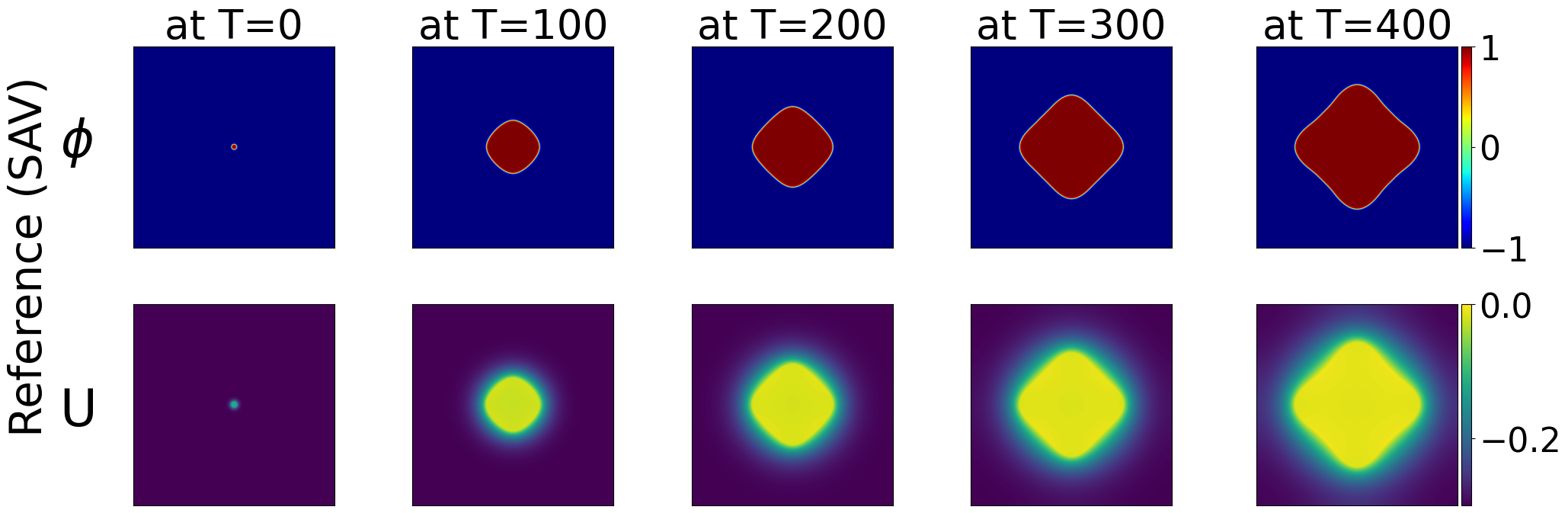}
  \includegraphics[width=.9\textwidth]{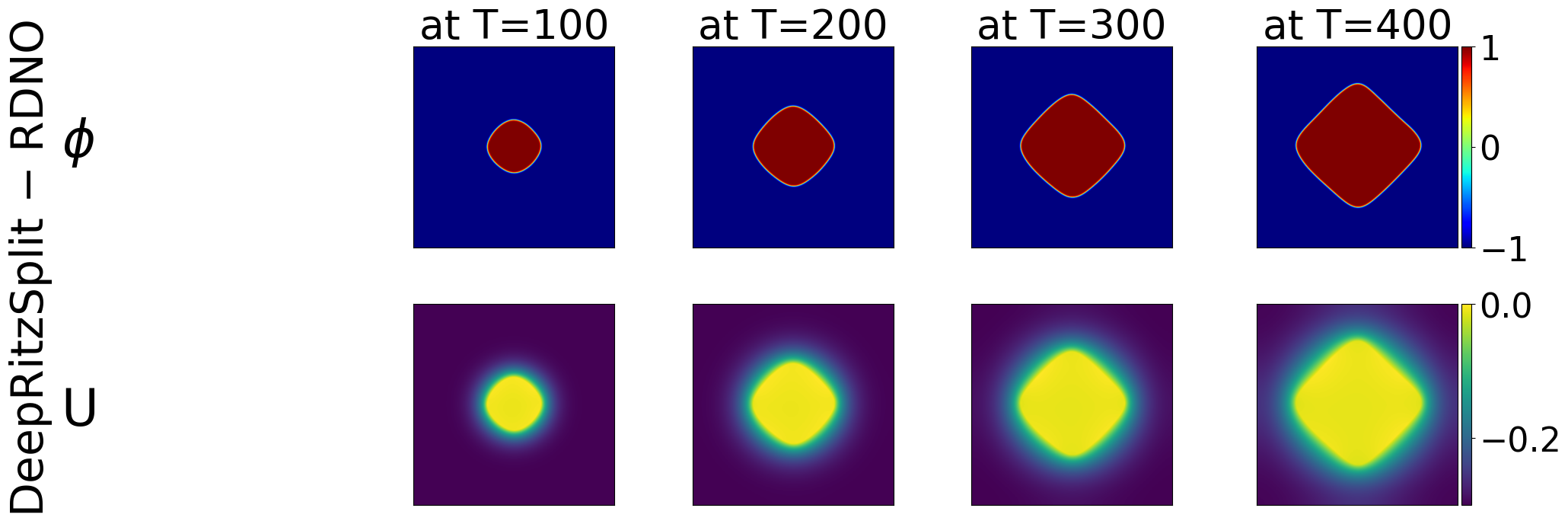}
  \includegraphics[width=.9\textwidth]{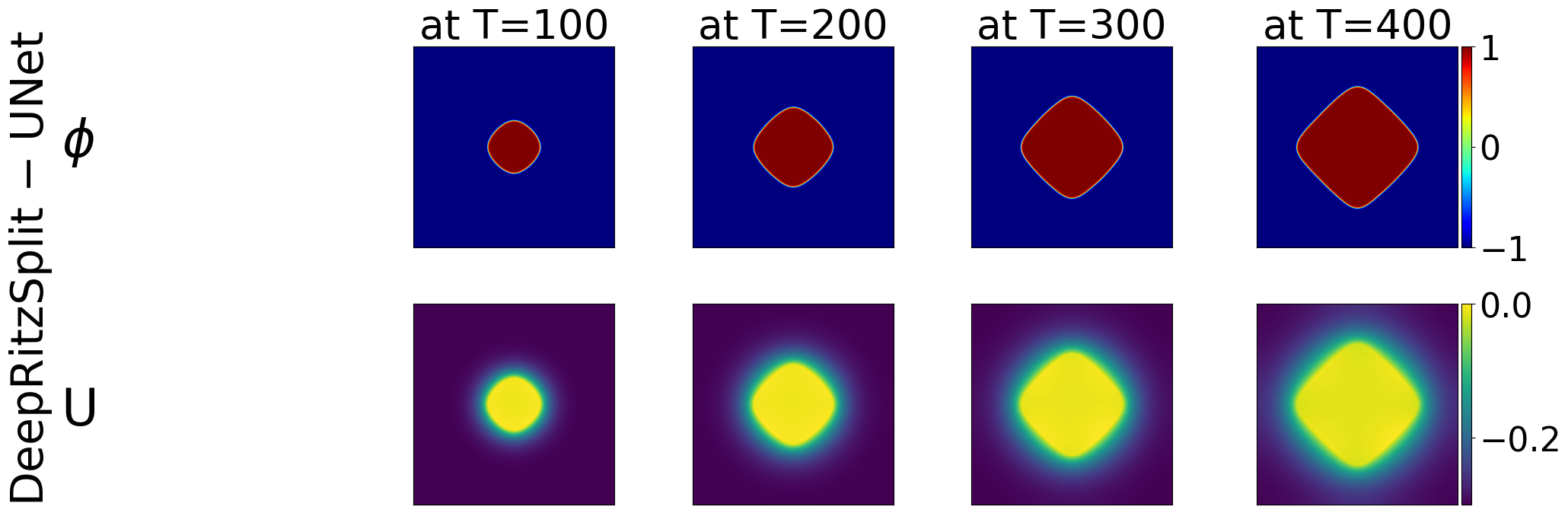}
  \includegraphics[width=.9\textwidth]{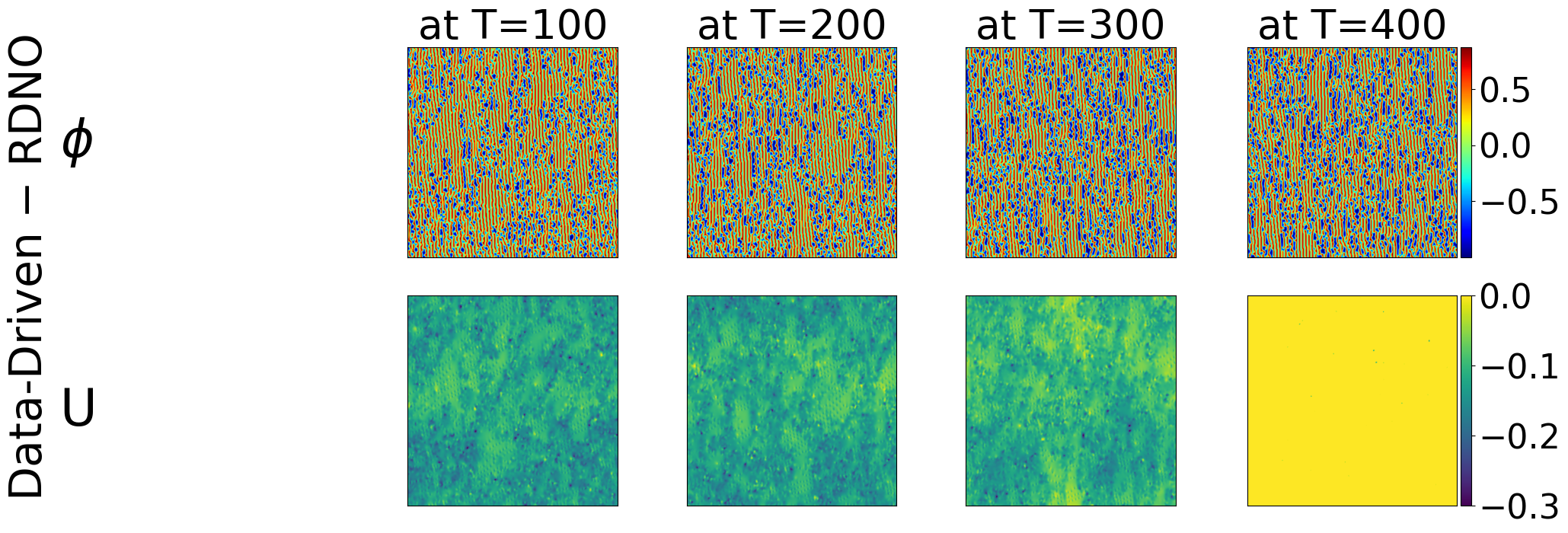}
  \caption{Dendritic growth predictions with anisotropic strength $\sigma=0.01$}
  \label{fig:SAV_RDNO_UNet-one-gem-eps_m001.png}
\end{figure}

\end{document}